\documentclass[12pt]{article}
\usepackage{amsfonts}
\usepackage{amsmath,amssymb,latexsym}
\usepackage{graphics}

\newtheorem{thrm}{Theorem}[section]

\newtheorem{dfnt}[thrm]{Definition}
\newtheorem{prot}[thrm]{Proposition}

\newtheorem{corl}[thrm]{Corollary}
\newtheorem{remk}[thrm]{Remark}
\newtheorem{exa}[thrm]{Example}

\numberwithin{equation}{section}

\begin{document}

\title{SOME UNIVERSAL NONLINEAR INEQUALITIES}
\author{Hamzeh Agahi\thanks{%
e-mail: h\_agahi@aut.ac.ir; h\_agahi@yahoo.com (H. Agahi)} \\
{\small \textit{Department of Statistics, Faculty of Mathematics and
Computer Science,}}\\
{\small \textit{Amirkabir University of Technology (Tehran Polytechnic),
424, Hafez Ave., Tehran 15914, Iran}}\\
}
\date{}
\maketitle

\begin{abstract}
In this paper, new versions of Chebyshev's, Minkowski's and H\"{o}lder's
type inequalities are studied by using a monotone measure-base universal
integral on an arbitrary measurable space. This paper generalizes some
previous results obtained by many researchers. \newline

\textit{Keywords:} Monotone measure; Universal integral; Chebyshev's
inequality; Minkowski's inequality; H\"{o}lder's inequality.
\end{abstract}

\section{Introduction}

Observe that in the last few years, there were introduced and discussed
several inequalities for non-classical integrals, thus developing a
theoretical background for further applications. Inequalities are at the
heart of the mathematical analysis of various problems in machine learning
and made it possible to derive new efficient algorithms.

In this paper, new versions of Chebyshev's, Minkowski's and H\"{o}lder's
type inequalities for universal integral on abstract spaces are studied in
rather general form, thus generalizing the results of \cite%
{Agahi2,Aga09,FloRom,MesOuy09,OuyFanWan08,OuyMesAg,OuyMes09}. Many nonlinear
systems are built by non-classical techniques, and thus we believe that our
results will prove their usefulness in flourishing areas, such as the
economy and decision making, among others.

The paper is organized as follows. In the next section, we briefly recall
some preliminaries and summarization of some previous known results. In
Section 3, we will focus on some interesting integral inequalities,
including Chebyshev's inequality, H\"{o}lder's inequality and Minkowski's
inequality for universal integral. Section 4 includes reverse previous
inequalities for semiconormed fuzzy integrals. Finally, a conclusion is
given.

\section{Universal integral}

\hspace{0.5cm}In this section, we are going to review some well-known known
results from universal integral. For the convenience of the reader, we
provide in this section a summary of the mathematical notations and
definitions used in this paper (see \textrm{\cite{KleMesPap08}}).

\begin{dfnt}
\label{th2-1}\textrm{\cite{KleMesPap08} }A monotone measure $m$ on a
measurable space $(X,\mathcal{A})$ is a function $m:\mathcal{A}\rightarrow
\lbrack 0,\infty ]$ satisfying \newline
(i) $m\left( \phi \right) =0,$ \newline
(ii) $m(X)>0$, \newline
(iii) $m(A)\leq m(B)$ whenever $A\subseteq B$.
\end{dfnt}

Note that a monotone measure is not necessarily $\sigma -$additive. This
concept goes back to M. Sugeno \cite{Sug} (where also the continuity of the
measures was required). To be precise, normed monotone measures on $(X,%
\mathcal{A})$, i.e., monotone measures satisfying $m(X)=1$, are also called
fuzzy measures \cite{ger,Sug,wang}, depending on the context.

For a fixed measurable space $(X,\mathcal{A})$, i.e., a non-empty set $X$
equipped with a $\sigma $-algebra $\mathcal{A}$, recall that a function $%
f:X\rightarrow \lbrack 0,\infty ]$ is called $\mathcal{A}$-measurable if,
for each $B$ $\in $ $\mathcal{B}\left( [0,\infty ]\right) $, the $\sigma $%
-algebra of Borel subsets of $[0,\infty ]$, the preimage $f^{-1}(B)$ is an
element of $\mathcal{A}$. We shall use the following notions:

\begin{dfnt}
\textrm{\cite{KleMesPap08} }Let $(X,\mathcal{A})$ be a measurable space. 
\newline
(i) $\mathcal{F}^{(X,\mathcal{A})}$ denotes the set of all $\mathcal{A}$%
-measurable functions $f:X\rightarrow \lbrack 0,\infty ]$; \newline
(ii) For each number $a\in $ $(0,\infty ]$, $\mathcal{M}_{a}^{(X,\mathcal{A}%
)}$ a denotes the set of all monotone measures (in the sense of Definition %
\ref{th2-1}) satisfying $m(X)=a$; and we take 
\begin{equation*}
\mathcal{M}^{(X,\mathcal{A})}=\bigcup_{a\in (0,\infty ]}\mathcal{M}_{a}^{(X,%
\mathcal{A})}.
\end{equation*}%
Let $\mathcal{S}$ be the class of all measurable spaces, and take 
\begin{equation*}
\mathcal{D}_{[0,\infty ]}=\bigcup_{\left( X,\mathcal{A}\right) \in \mathcal{S%
}}\mathcal{M}^{(X,\mathcal{A})}\times \mathcal{F}^{(X,\mathcal{A})}.
\end{equation*}
\end{dfnt}

The Choquet \textrm{\cite{ch}}, Sugeno \textrm{\cite{Sug} }and Shilkret 
\textrm{\cite{shil}} integrals (see also \textrm{\cite{BenMesViv,Pap02}}),
respectively, are given, for any measurable space $\left( X,\mathcal{A}%
\right) $, for any measurable function $f\in $ $\mathcal{F}^{(X,\mathcal{A}%
)} $ and for any monotone measure $m$ $\in \mathcal{M}^{(X,\mathcal{A})}$,
i.e., for any $\left( m,f\right) \in \mathcal{D}_{[0,\infty ]},$ by 
\begin{eqnarray}
\mathbf{Su}(m,f) &=&\sup \left\{ \min \left( t,m\left( \left\{ f\geq
t\right\} \right) \right) \;|\;t\in (0,\infty ])\right\} ,  \label{s} \\
\mathbf{Sh}(m,f) &=&\sup \left\{ t.m\left( \left\{ f\geq t\right\} \right)
\;|\;t\in (0,\infty ])\right\} ,  \label{t}
\end{eqnarray}%
where the convention $0.\infty =0$ is used. All these integrals map $%
\mathcal{M}^{(X,\mathcal{A})}$ $\times \mathcal{F}^{(X,\mathcal{A})}$ into $%
[0,\infty ]$ independently of $(X,\mathcal{A}).$ We remark that fixing an
arbitrary $m\in \mathcal{M}^{(X,\mathcal{A})}$, they are non-decreasing
functions from $\mathcal{F}^{(X,\mathcal{A})}$ into $[0,\infty ]$, and
fixing an arbitrary $f\in \mathcal{F}^{(X,\mathcal{A})}$, they are
non-decreasing functions from $\mathcal{M}^{(X,\mathcal{A})}$ into $%
[0,\infty ]$.

We stress the following important common property for all three integrals
from (\ref{s}) and (\ref{t}). Namely, these integrals does not make
difference between the pairs $\left( m_{1},f_{1}\right) ,\left(
m_{2},f_{2}\right) $ $\in \mathcal{D}_{[0,\infty ]}$ which satisfy, for all
for all $t\in (0,\infty ],$%
\begin{equation*}
m_{1}(\left\{ f_{1}\geq t\right\} )=m_{2}(\left\{ f_{2}\geq t\right\} ).
\end{equation*}%
Therefore, such equivalence relation between pairs of measures and functions
was introduced in \textrm{\cite{KleMesPap08}}.

\begin{dfnt}
\textrm{\ }Two pairs $(m_{1},f_{1})\in \mathcal{M}^{(X_{1},\mathcal{A}_{1})}$
$\times \mathcal{F}^{(X_{1},\mathcal{A}_{1})}$ and $(m_{2},f_{2})\in 
\mathcal{M}^{(X_{2},\mathcal{A}_{2})}$ $\times \mathcal{F}^{(X_{2},\mathcal{A%
}_{2})}$ satisfying 
\begin{equation*}
m_{1}(\left\{ f_{1}\geq t\right\} )=m_{2}(\left\{ f_{2}\geq t\right\} )\text{
for all }t\in (0,\infty ],
\end{equation*}%
will be called integral equivalent, in symbols 
\begin{equation*}
\left( m_{1},f_{1}\right) \sim \left( m_{2},f_{2}\right) .
\end{equation*}
\end{dfnt}

To introduce the notion of the universal integral we shall need instead of
the usual plus and product more general real operations.

\begin{dfnt}
\textrm{\cite{SugMur} }A function $\otimes \colon \lbrack 0,\infty ]\sp{2}%
\rightarrow \lbrack 0,\infty ]$ is called a pseudo-multiplication if it
satisfies the following properties:\newline
(i) it is non-decreasing in each component, i.e., for all $%
a_{1},a_{2},b_{1},b_{2}\in \lbrack 0,\infty ]$ with $a_{1}\leq a_{2}$ and $%
b_{1}\leq b_{2}$ we have $a_{1}\otimes b_{1}\leq a_{2}\otimes b_{2}$;\newline
(ii) $0$ is an annihilator of , i.e., for all $a\in \lbrack 0,\infty ]$ we
have $a\otimes 0=0\otimes a=0$;\newline
(iii) has a neutral element different from $0$, i.e., there exists an $e\in
(0,\infty ]$ such that, for all $a\in \lbrack 0,\infty ]$, we have $a\otimes
e=e\otimes a=a.$
\end{dfnt}

There is neither a smallest nor a greatest pseudo-multiplication on $%
[0,\infty ]$. But, if we fix the neutral element $e\in (0,\infty ]$, then
the smallest pseudo-multiplication $\otimes _{e}$ and the greatest
pseudo-multiplication $\otimes ^{e}$ with neutral element $e$ are given by 
\begin{equation*}
a\otimes _{e}b=\left\{ 
\begin{array}{cc}
0 & if\;\left( a,b\right) \in \lbrack 0,e)^{2}, \\ 
\max \left( a,b\right) & if\;\left( a,b\right) \in \lbrack e,\infty ]^{2},
\\ 
\min \left( a,b\right) & otherwise,%
\end{array}%
\right.
\end{equation*}%
and 
\begin{equation*}
a\otimes ^{e}b=\left\{ 
\begin{array}{cc}
\min \left( a,b\right) & if\;\min \left( a,b\right) =0\text{ or }\left(
a,b\right) \in (0,e]^{2}, \\ 
\infty & if\;\left( a,b\right) \in (e,\infty ]^{2}, \\ 
\max \left( a,b\right) & otherwise.%
\end{array}%
\right.
\end{equation*}

Restricting to the interval $[0,1]$ a pseudo-multiplication and a
pseudo-addition with additional properties of associativity and
commutativity can be considered as the $t$-norm $T$ and the $t$-conorms $S$
(see \textrm{\cite{KleMesPap}}), respectively.

For a given pseudo-multiplication on $[0,\infty ]$, we suppose the existence
of a pseudo-addition $\oplus \colon \lbrack 0,\infty ]\sp{2}\rightarrow
\lbrack 0,\infty ]$ which is continuous, associative, non-decreasing and has 
$0$ as neutral element (then the commutativity of follows, see \textrm{\cite%
{KleMesPap}}), and which is left-distributive with respect to $\otimes $
i.e., for all $a,b,c$ $\in \lbrack 0,\infty ]$ we have $(a\oplus b)\otimes
c=(a\oplus c)\otimes (b\oplus c).$ The pair $(\oplus ,\otimes )$ is then
called an integral operation pair, see \textrm{\cite{BenMesViv,KleMesPap08}.}

Each of the integrals mentioned in (\ref{s}) and (\ref{t}) maps $\mathcal{D}%
_{[0,\infty ]}$ into $[0,\infty ]$ and their main properties can be covered
by the following common integral given in \textrm{\cite{KleMesPap08}}.

\begin{dfnt}
A function $\mathbf{I}\colon \mathcal{D}_{[0,\infty ]}\rightarrow \lbrack
0,\infty ]$ is called a universal integral if the following axioms hold:%
\newline
(I1) For any measurable space $(X,\mathcal{A})$, the restriction of the
function $\mathbf{I}$ to $\mathcal{M}^{(X,\mathcal{A})}$ $\times \mathcal{F}%
^{(X,\mathcal{A})}$ is non-decreasing in each coordinate;\newline
(I2) there exists a pseudo-multiplication $\otimes \colon \lbrack 0,\infty ]%
\sp{2}\rightarrow \lbrack 0,\infty ]$ such that for all pairs $(m,c.\mathbf{1%
}_{A})$ $\in \mathcal{D}_{[0,\infty ]}$ 
\begin{equation*}
\mathbf{I}(m,c.\mathbf{1}_{A})=c\otimes m(A);
\end{equation*}%
(I3) for all integral equivalent pairs $\left( m_{1},f_{1}\right) ,\left(
m_{2},f_{2}\right) $ $\in \mathcal{D}_{[0,\infty ]}$ we have $\mathbf{I}%
\left( m_{1},f_{1}\right) =\mathbf{I}\left( m_{2},f_{2}\right) .$
\end{dfnt}

By Proposition 3.1 from \textrm{\cite{KleMesPap08}} we have the following
important characterization.

\begin{thrm}
Let $\otimes \colon \lbrack 0,\infty ]\sp{2}\rightarrow \lbrack 0,\infty ]$
be a pseudo-multiplication on $[0,\infty ]$. Then the smallest universal
integral $\mathbf{I}$ based on $\otimes $ is given by 
\begin{equation*}
\mathbf{I}_{\otimes }\left( m,f\right) =\sup \left\{ t\otimes m\left(
\left\{ f\geq t\right\} \right) \;|\;t\in (0,\infty ])\right\} .
\end{equation*}
\end{thrm}

Specially, we have $\mathbf{Su}=$ $\mathbf{I}_{Min}$ and $\mathbf{Sh}=$ $%
\mathbf{I}_{Prod}$, where the pseudo-multiplications $Min$ and $Prod$ are
given (as usual) by $Min(a,b)=min(a,b)$ and $Prod(a,b)=a.b$. Note that the
nonlinearity of the Sugeno integral $\mathbf{Su}$ (see, e.g., \textrm{\cite%
{Kle,Mes2008}}) implies that universal integrals are also nonlinear, in
general.

\begin{prot}
There exists the smallest universal integral $\mathbf{I}_{\otimes _{e}}$
among all universal integrals satisfying the conditions \newline
(i) for each $m\in \mathcal{M}_{e}^{(X,\mathcal{A})}$ and each $c\in \lbrack
0,\infty ]$ we have $\mathbf{I}(m,c.\mathbf{1}_{X})=c$,\newline
(ii) for each $m\in \mathcal{M}^{(X,\mathcal{A})}$ and each $A\in \mathcal{A}
$ we have $\mathbf{I}(m,e.\mathbf{1}_{X})=m(A)$, given by%
\begin{equation*}
\mathbf{I}_{\otimes _{e}}\left( m,f\right) =\max \left\{ m\left( \left\{
f\geq e\right\} \right) ,essinf_{m}f\right\}
\end{equation*}%
where $essinf_{m}f=\sup \left\{ t\in \lbrack 0,\infty ]\;|\;m\left( \left\{
f\geq t\right\} \right) =m\left( X\right) \right\} $.
\end{prot}

Restricting now to the unit interval $[0,1]$ we shall consider functions $%
f\in \mathcal{F}^{(X,\mathcal{A})}$ satisfying $Ran(f)\subseteq \lbrack 0,1]$
(in which case we shall write shortly $f\in \mathcal{F}_{\left[ 0,1\right]
}^{(X,\mathcal{A})}$). Observe that, in this case, we have the restriction
of the pseudo-multiplication $\otimes $ to $[0,1]^{2}$ (called a semicopula
or a conjunctor, i.e., a binary operation $\circledast \colon \lbrack 0,1]%
\sp{2}\rightarrow \lbrack 0,1]$ which is non-decreasing in both components,
has $1$ as neutral element and satisfies $a\circledast b\leq $ $\min (a,b)$
for all $(a,b)\in \lbrack 0,1]^{2},$ see \cite{Bas,DurSem}), and universal
integrals are restricted to the class $\mathcal{D}_{[0,1]}=\bigcup_{\left( X,%
\mathcal{A}\right) \in \mathcal{S}}\mathcal{M}^{(X,\mathcal{A})}\times 
\mathcal{F}_{\left[ 0,1\right] }^{(X,\mathcal{A})}$. In a special case, for
a fixed strict $t$-norm $T$, the corresponding universal integral $\mathbf{I}%
_{T}$ is the so-called Sugeno-Weber integral \cite{Web}. The smallest
universal integral $\mathbf{I}_{\circledast }$ on the $[0,1]$ scale related
to the semicopula $\circledast $ is given by

\begin{equation*}
\mathbf{I}_{\circledast }\left( m,f\right) =\sup \left\{ t\circledast
m\left( \left\{ f\geq t\right\} \right) \;|\;t\in \lbrack 0,1])\right\} .
\end{equation*}%
This type of integral was called seminormed integral in \cite{SurGil}.%
\newline
Before starting our main results we need the following definitions:

\begin{dfnt}
Functions $f,g\colon X\rightarrow \mathbb{R}$ are said to be comonotone if
for all $x,y\in X$, 
\begin{equation*}
(f(x)-f(y))(g(x)-g(y))\geq 0,
\end{equation*}%
and $f$ and $g$ are said to be countermonotone if for all $x,y\in X$, 
\begin{equation*}
(f(x)-f(y))(g(x)-g(y))\leq 0.
\end{equation*}
\end{dfnt}

The comonotonicity of functions $f$ and $g$ is equivalent to the
nonexistence of points $x,y\in X$ such that $f(x)<f(y)$ and $g(x)>g(y)$.
Similarly, if $f$ and $g$ are countermonotone then $f(x)<f(y)$ and $%
g(x)<g(y) $ cannot happen. Observe that the concept of comonotonicity was
first introduced in \textrm{\cite{De}}.\newline

\begin{dfnt}
Let $A,B\colon \lbrack 0,\infty ]\sp{2}\rightarrow \lbrack 0,\infty ]$ be
two binary operations. Recall that $A$ dominates $B$ (or $B$ is dominated by 
$A$), denoted by $A\gg B$, if 
\begin{equation*}
A(B(a,b),B(c,d))\geq B(A(a,c),A(b,d))
\end{equation*}%
holds for any $a,b,c,d\in \lbrack 0,\infty ]$.
\end{dfnt}

\begin{dfnt}
Let $\star \colon \lbrack 0,\infty ]\sp{2}\rightarrow \lbrack 0,\infty ]$ be
a binary operation and consider $\varphi :[0,\infty ]\rightarrow \lbrack
0,\infty ]$. Then we say that $\varphi $ is subdistributive over $\star $ if 
\begin{equation*}
\varphi (x\star y)\leq \varphi (x)\star \varphi (y)
\end{equation*}%
for all $x,y\in \lbrack 0,\infty ]$. Analogously, we say that $\varphi $ is
superdistributive over $\star $ if 
\begin{equation*}
\varphi (x\star y)\geq \varphi (x)\star \varphi (y)
\end{equation*}%
for all $x,y\in \lbrack 0,\infty ]$.
\end{dfnt}

\section{On some advanced type inequalities for universal integral}

\hspace{5mm} Now, we state the main result of this paper.

\begin{thrm}
\label{thfi} Let a non-decreasing $n$-place function $H:[0,\infty
)^{n}\rightarrow \lbrack 0,\infty )$ such that $H$ be continuous. If $\
\otimes \colon \lbrack 0,\infty ]\sp{n}\rightarrow \lbrack 0,\infty ]$ is
the pseudo-multiplication with neutral element $e\in (0,\infty ]$, satisfies 
\begin{eqnarray*}
&&U_{0}^{-1}\left[ U_{0}\left( H\left( \psi _{1}\left( a_{1}\right) ,\psi
_{2}\left( a_{2}\right) ,...,\psi _{n}\left( a_{n}\right) \right) \right)
\otimes c\right] \\
&&\overset{}{\geq }\left[ 
\begin{array}{c}
H\left( \psi _{1}\left( U_{1}^{-1}\left[ \left( U_{1}\left( a_{1}\right)
\right) \otimes c\right] \right) ,\psi _{2}\left( a_{2}\right) ,...,\psi
_{n}\left( a_{n}\right) \right) \\ 
\vee H\left( \psi _{1}\left( a_{1}\right) ,\psi _{2}\left( U_{2}^{-1}\left[
\left( U_{2}\left( a_{2}\right) \right) \otimes c\right] \right) ,\psi
_{2}\left( a_{3}\right) ,...,\psi _{n}\left( a_{n}\right) \right) \\ 
\vee ...\vee H\left( \psi _{1}\left( a_{1}\right) ,\psi _{2}\left(
a_{2}\right) ,...,\psi _{n-1}\left( a_{n-1}\right) ,\psi _{n}\left(
U_{n}^{-1}\left[ \left( U_{n}\left( a_{n}\right) \right) \otimes c\right]
\right) \right) ,%
\end{array}%
\right]
\end{eqnarray*}%
then for any system $U_{0},U_{1},...,U_{n}:[0,\infty )\rightarrow \lbrack
0,\infty )$ of continuous strictly increasing functions, and any system $%
\psi _{1},\psi _{2},...,\psi _{n}:[0,\infty )\rightarrow \lbrack 0,\infty )$
of continuous increasing functions and any comontone system $%
f_{1},f_{2},...,f_{n}\in $ $\mathcal{F}^{(X,\mathcal{A})}$ and a monotone
measure $m$ $\in \mathcal{M}^{(X,\mathcal{A})}$ such that $b\otimes m\left(
X\right) \leq b$ for all $b\in \left[ 0,\infty \right] $ and $\mathbf{I}%
_{\otimes }\left( m,U_{i}\left( f_{i}\right) \right) <\infty $ for all $%
i=1,2,...n$, it holds%
\begin{equation*}
U_{0}^{-1}[\mathbf{I}_{\otimes }\left( m,U_{0}[H\left( \psi _{1}\left(
f_{1}\right) ,...,\psi _{n}\left( f_{n}\right) \right) ]\right) ]\overset{}{%
\geq }H\left[ \psi _{1}\left( U_{1}^{-1}\left( \mathbf{I}_{\otimes }\left(
m,U_{1}\left( f_{1}\right) \right) \right) \right) ,...,\psi _{n}\left(
U_{n}^{-1}\left( \mathbf{I}_{\otimes }\left( m,U_{n}\left( f_{n}\right)
\right) \right) \right) \right] .
\end{equation*}
\end{thrm}

\noindent \textbf{Proof.}\ Let $e\in (0,\infty ]$ be the neutral element of $%
\otimes $ and $\mathbf{I}_{\otimes }\left( m,U_{i}\left( f_{i}\right)
\right) =p_{i}<\infty $ for all $i=1,2,...,n$. So, for any $\varepsilon >0$,
there exist $p_{i(\varepsilon )}$ such that 
\begin{equation*}
m(\left\{ U_{i}\left( f_{i}\right) \geq p_{i(\varepsilon )}\right\}
)=m(\left\{ f_{i}\geq U_{i}^{-1}\left( p_{i(\varepsilon )}\right) \right\}
)=M\sb{i},
\end{equation*}%
where $p_{i(\varepsilon )}\otimes M\sb{i}\geq p_{i}-\varepsilon $ for all $%
i=1,2,...,n.$ Then,%
\begin{equation*}
\psi _{i}\left( U_{i}^{-1}\left[ p_{i\left( \varepsilon \right) }\otimes M%
\sb{i}\right] \right) \geq \psi _{i}\left( U_{i}^{-1}\left[
p_{i}-\varepsilon \right] \right) ,\text{ for all }i=1,2,...,n.
\end{equation*}%
Then,%
\begin{equation*}
\psi _{i}\left( U_{i}^{-1}\left[ p_{i\left( \varepsilon \right) }\right]
\right) \geq \psi _{i}\left( U_{i}^{-1}\left[ p_{i\left( \varepsilon \right)
}\otimes m\left( X\right) \right] \right) \geq \psi _{i}\left( U_{i}^{-1}%
\left[ p_{i}-\varepsilon \right] \right) ,\text{ for all }i=1,2,...,n.
\end{equation*}%
The comonotonicity of $f_{1},f_{2},...,f_{n}$ and the monotonicity of $H$
imply that 
\begin{eqnarray*}
&&m\left( \{U_{0}\left( H\left( \psi _{1}\left( f_{1}\right) ,...,\psi
_{n}\left( f_{n}\right) \right) \right) \geq U_{0}\left( H\left( \psi
_{1}\left( U_{1}^{-1}\left( p_{1(\varepsilon )}\right) \right) ,...,\psi
_{n}\left( U_{n}^{-1}\left( p_{n(\varepsilon )}\right) \right) \right)
\right) \}\right) \\
&=&m(\{H\left( \psi _{1}\left( f_{1}\right) ,...,\psi _{n}\left(
f_{n}\right) \right) \geq H\left( \psi _{1}\left( U_{1}^{-1}\left(
p_{1(\varepsilon )}\right) \right) ,...,\psi _{n}\left( U_{n}^{-1}\left(
p_{n(\varepsilon )}\right) \right) \right) \}) \\
&\geq &m(\left\{ f_{1}\geq U_{1}^{-1}\left( p_{1(\varepsilon )}\right)
\right\} )\wedge m(\left\{ f_{2}\geq U_{2}^{-1}\left( p_{2(\varepsilon
)}\right) \right\} )\wedge ....\wedge m(\left\{ f_{n}\geq U_{n}^{-1}\left(
p_{n(\varepsilon )}\right) \right\} ) \\
&=&M_{1}\wedge M_{2}\wedge ...\wedge M_{n}.
\end{eqnarray*}%
Hence 
\begin{eqnarray*}
&&U_{0}^{-1}\left[ \sup \left( t\otimes m(\{U_{0}\left( H\left( \psi
_{1}\left( f_{1}\right) ,...,\psi _{n}\left( f_{n}\right) \right) \right)
\geq t\})\;|\;t\in (0,\infty ])\right) \right] \\
&&\overset{}{\geq }U_{0}^{-1}\left( \left[ 
\begin{array}{c}
U_{0}\left( H\left( \psi _{1}\left( U_{1}^{-1}\left( p_{1(\varepsilon
)}\right) \right) ,...,\psi _{n}\left( U_{n}^{-1}\left( p_{n(\varepsilon
)}\right) \right) \right) \right) \otimes \\ 
m(\{U_{0}\left( H\left( \psi _{1}\left( f_{1}\right) ,...,\psi _{n}\left(
f_{n}\right) \right) \right) \geq U_{0}\left( H\left( \psi _{1}\left(
U_{1}^{-1}\left( p_{1(\varepsilon )}\right) \right) ,...,\psi _{n}\left(
U_{n}^{-1}\left( p_{n(\varepsilon )}\right) \right) \right) \right) \})%
\end{array}%
\right] \right) \\
&&\overset{}{\geq }U_{0}^{-1}\left( \left[ U_{0}\left( H\left( \psi
_{1}\left( U_{1}^{-1}\left( p_{1(\varepsilon )}\right) \right) ,...,\psi
_{n}\left( U_{n}^{-1}\left( p_{n(\varepsilon )}\right) \right) \right)
\right) \otimes \left( M_{1}\wedge M_{2}\wedge ...\wedge M_{n}\right) \right]
\right) \\
&&\overset{}{=}\left( 
\begin{array}{c}
U_{0}^{-1}\left[ U_{0}\left( H\left( \psi _{1}\left( U_{1}^{-1}\left(
p_{1(\varepsilon )}\right) \right) ,...,\psi _{n}\left( U_{n}^{-1}\left(
p_{n(\varepsilon )}\right) \right) \right) \right) \otimes M_{1}\right] \\ 
\wedge U_{0}^{-1}\left[ U_{0}\left( H\left( \psi _{n}\left( U_{1}^{-1}\left(
p_{1(\varepsilon )}\right) \right) ,...,\psi _{n}\left( U_{n}^{-1}\left(
p_{n(\varepsilon )}\right) \right) \right) \right) \otimes M_{2}\right] \\ 
\wedge ...\wedge U_{0}^{-1}\left[ U_{0}\left( H\left( \psi _{1}\left(
U_{1}^{-1}\left( p_{1(\varepsilon )}\right) \right) ,...,\psi _{n}\left(
U_{n}^{-1}\left( p_{n(\varepsilon )}\right) \right) \right) \right) \otimes
M_{n}\right]%
\end{array}%
\right) \\
&&\overset{}{\geq }\left( 
\begin{array}{c}
H\left( \psi _{1}\left( U_{1}^{-1}\left[ p_{1\left( \varepsilon \right)
}\otimes M_{1}\right] \right) ,\psi _{2}\left( U_{2}^{-1}\left[ p_{2\left(
\varepsilon \right) }\right] \right) ,...,\psi _{n}\left( U_{n}^{-1}\left[
p_{n\left( \varepsilon \right) }\right] \right) \right) \\ 
\wedge H\left( \psi _{1}\left( U_{1}^{-1}\left[ p_{1\left( \varepsilon
\right) }\right] \right) ,\psi _{2}\left( U_{2}^{-1}\left[ p_{2\left(
\varepsilon \right) }\otimes M_{2}\right] \right) ,...,\psi _{n}\left(
U_{n}^{-1}\left[ p_{n\left( \varepsilon \right) }\right] \right) \right) \\ 
\wedge ...\wedge H\left( \psi _{1}\left( U_{1}^{-1}\left[ p_{1\left(
\varepsilon \right) }\right] \right) ,...,\psi _{n-1}\left( U_{n-1}^{-1}%
\left[ p_{\left( n-1\right) \left( \varepsilon \right) }\right] \right)
,\psi _{n}\left( U_{n}^{-1}\left[ p_{n\left( \varepsilon \right) }\otimes
M_{n}\right] \right) \right)%
\end{array}%
\right) \\
&&\overset{}{\geq }\left( 
\begin{array}{c}
H\left( \psi _{1}\left( U_{1}^{-1}\left[ p_{1}-\varepsilon \right] \right)
,\psi _{2}\left( U_{2}^{-1}\left[ p_{2\left( \varepsilon \right) }\right]
\right) ,...,\psi _{n}\left( U_{n}^{-1}\left[ p_{n\left( \varepsilon \right)
}\right] \right) \right) \\ 
\wedge H\left( \psi _{1}\left( U_{1}^{-1}\left[ p_{1\left( \varepsilon
\right) }\right] \right) ,\psi _{2}\left( U_{2}^{-1}\left[ p_{2}-\varepsilon %
\right] \right) ,...,\psi _{n}\left( U_{n}^{-1}\left[ p_{n\left( \varepsilon
\right) }\right] \right) \right) \\ 
\wedge ...\wedge H\left( \psi _{1}\left( U_{1}^{-1}\left[ p_{1\left(
\varepsilon \right) }\right] \right) ,...,\psi _{n-1}\left( U_{n-1}^{-1}%
\left[ p_{\left( n-1\right) \left( \varepsilon \right) }\right] \right)
,\psi _{n}\left( U_{n}^{-1}\left[ p_{n}-\varepsilon \right] \right) \right)%
\end{array}%
\right) \\
&&\overset{}{\geq }H\left( \psi _{1}\left( U_{1}^{-1}\left[
p_{1}-\varepsilon \right] \right) ,\psi _{2}\left( U_{2}^{-1}\left[
p_{2}-\varepsilon \right] \right) ,...,\psi _{n}\left( U_{n}^{-1}\left[
p_{n}-\varepsilon \right] \right) \right) ,
\end{eqnarray*}%
whence $U_{0}^{-1}[\mathbf{I}_{\otimes }\left( m,U_{0}[H\left( \psi
_{1}\left( f_{1}\right) ,...,\psi _{n}\left( f_{n}\right) \right) ]\right)
]\geq H\left( \psi _{1}\left( U_{1}^{-1}\left[ p_{1}\right] \right) ,\psi
_{2}\left( U_{2}^{-1}\left[ p_{2}\right] \right) ,...,\psi _{n}\left(
U_{n}^{-1}\left[ p_{n}\right] \right) \right) $ follows from the continuity
of $H,\psi _{i},U_{i}$ for all $i,$ and the arbitrariness of $\varepsilon $.
And the theorem is proved. $\Box $

\begin{remk}
\label{remark2} (i) If $m(X)=e$, then the condition $b\otimes m\left(
X\right) \leq b$ for all $b\in \lbrack 0,\infty ]$ holds readily.\newline
(ii) We can replace the condition \textquotedblleft $b\otimes m\left(
X\right) \leq b$ for all $b\in \lbrack 0,\infty ]"$ with \textquotedblleft $%
b\otimes a\leq b$ for all $a,b\in \lbrack 0,\infty ]".$
\end{remk}

\begin{corl}
\label{cor3.1} Let $f,g\in $ $\mathcal{F}^{(X,\mathcal{A})}$ be two
comonotone measurable functions and $\otimes \colon \lbrack 0,\infty ]\sp{2}%
\rightarrow \lbrack 0,\infty ]$ be the pseudo-multiplication with neutral
element $e\in (0,\infty ]$ and $m$ $\in \mathcal{M}^{(X,\mathcal{A})}$ be a
monotone measure such that $a\otimes m\left( X\right) \leq a$ for all $a\in %
\left[ 0,\infty \right] ,$ $\mathbf{I}_{\otimes }\left( m,U_{1}\left(
f\right) \right) $ and $\mathbf{I}_{\otimes }\left( m,U_{2}\left( g\right)
\right) $ are finite and $U_{i}:[0,\infty )\rightarrow \lbrack 0,\infty ),$ $%
i=0,1,2$ be continuous strictly increasing functions. Let $\star \colon
\lbrack 0,\infty )\sp{2}\rightarrow \lbrack 0,\infty )$ be continuous and
nondecreasing in both arguments and $\psi :[0,\infty )\rightarrow \lbrack
0,\infty )$ be continuous and strictly increasing function. If $\ $ 
\begin{eqnarray}
U_{0}^{-1}\left[ U_{0}\left( \psi \left( a\right) \star \psi \left( b\right)
\right) \otimes c\right] &\geq &\left[ \psi \left( U_{1}^{-1}\left[ \left(
U_{1}\left( a\right) \right) \otimes c\right] \right) \star \psi \left(
b\right) \right]  \label{rtr} \\
&&\vee \left[ \psi \left( a\right) \star \psi \left( U_{2}^{-1}\left[ \left(
U_{2}\left( b\right) \right) \otimes c\right] \right) \right] ,  \notag
\end{eqnarray}%
then the inequality 
\begin{equation*}
U_{0}^{-1}[\mathbf{I}_{\otimes }\left( m,U_{0}[\left( \psi \left( f\right)
\star \psi \left( g\right) \right) ]\right) ]\geq \psi \left(
U_{1}^{-1}\left( \mathbf{I}_{\otimes }\left( m,U_{1}\left( f\right) \right)
\right) \right) \star \psi \left( U_{2}^{-1}\left( \mathbf{I}_{\otimes
}\left( m,U_{2}\left( g\right) \right) \right) \right)
\end{equation*}%
holds.
\end{corl}

Let $U_{i}\left( x\right) =\varphi _{i}\left( x\right) $ for all $i=0,1,2$
and $\psi \left( x\right) =x$ in Corollary \ref{cor3.1}. Then we have the
following result.

\begin{corl}
\label{cor3.2}Let $f,g\in $ $\mathcal{F}^{(X,\mathcal{A})}$ be two
comonotone measurable functions and $\otimes \colon \lbrack 0,\infty ]\sp{2}%
\rightarrow \lbrack 0,\infty ]$ be the pseudo-multiplication with neutral
element $e\in (0,\infty ]$ and $m$ $\in \mathcal{M}^{(X,\mathcal{A})}$ be a
monotone measure such that $a\otimes m\left( X\right) \leq a$ for all $a\in %
\left[ 0,\infty \right] ,$ $\mathbf{I}_{\otimes }\left( m,\varphi _{1}\left(
f\right) \right) $ and $\mathbf{I}_{\otimes }\left( m,\varphi _{2}\left(
g\right) \right) $ are finite. Let $\star \colon \lbrack 0,\infty )\sp{2}%
\rightarrow \lbrack 0,\infty )$ be continuous and nondecreasing in both
arguments and $\varphi _{i}:[0,\infty )\rightarrow \lbrack 0,\infty )$ $%
i=0,1,2$ be continuous strictly increasing functions. If $\ $ 
\begin{equation*}
\varphi _{0}^{-1}\left[ \varphi _{0}\left( p_{1}\star p_{2}\right) \otimes c%
\right] \geq \left[ \varphi _{1}^{-1}\left[ \left( \varphi _{1}\left(
p_{1}\right) \right) \otimes c\right] \star p_{2}\right] \vee \left[
p_{1}\star \varphi _{2}^{-1}\left[ \varphi _{2}\left( p_{2}\right) \otimes c%
\right] \right] ,
\end{equation*}%
then the inequality 
\begin{equation*}
\varphi _{0}^{-1}[\mathbf{I}_{\otimes }\left( m,\varphi _{0}\left( f\star
g\right) \right) ]\geq \varphi _{1}^{-1}\left( \mathbf{I}_{\otimes }\left(
m,\varphi _{1}\left( f\right) \right) \right) \star \varphi _{2}^{-1}\left( 
\mathbf{I}_{\otimes }\left( m,\varphi _{2}\left( g\right) \right) \right)
\end{equation*}%
holds.
\end{corl}

In an analogous way as in the proof of Theorem \ref{thfi} we have the
following result.

\begin{thrm}
\label{th3-2} Let $H:[0,\infty )^{n}\rightarrow \lbrack 0,\infty )$ be a
continuous and nondecreasing $n$-place function. If $\ \otimes \colon
\lbrack 0,\infty ]\sp{n}\rightarrow \lbrack 0,\infty ]$ is the
pseudo-multiplication on $[0,\infty ]$ with neutral element $e\in (0,\infty
] $ such that $a\otimes m\left( X\right) \leq a$ for all $a\in \lbrack
0,\infty ]$, satisfies 
\begin{eqnarray}
&&\left[ \left( H\left( p_{1},p_{2},...,p_{n}\right) \right) ^{\xi
_{0}}\otimes c\right] ^{\omega _{0}}\overset{}{\geq }H\left( \left(
p_{1}^{\xi _{1}}\otimes c\right) ^{\omega _{1}},p_{2},...,p_{n}\right) \vee
\label{tyr} \\
&&H\left( p_{1},\left( p_{2}^{\xi _{2}}\otimes c\right) ^{\omega
_{2}},p_{3},...,p_{n}\right) \vee ...\vee H\left(
p_{1},p_{2},...,p_{n-1},\left( p_{n}^{\xi _{n}}\otimes c\right) ^{\omega
_{n}}\right) ,  \notag
\end{eqnarray}%
then for any comontone system $f_{1},f_{2},...,f_{n}\in $ $\mathcal{F}^{(X,%
\mathcal{A})}$ and a monotone measure $m$ $\in \mathcal{M}^{(X,\mathcal{A})}$
such that $\mathbf{I}_{\otimes }\left( m,f_{i}^{\xi _{i}}\right) <\infty $
and $x^{\frac{1}{\xi _{i}\omega _{i}}}\geq x$ for all $x\in \lbrack 0,\infty
)$ and $i=1,2,...n$, it holds%
\begin{equation}
\left[ \mathbf{I}_{\otimes }\left( m,\left( H\left( f_{1},...,f_{n}\right)
\right) ^{\xi _{0}}\right) \right] ^{^{\omega _{0}}}\overset{}{\geq }H\left[
\left( \mathbf{I}_{\otimes }\left( m,f_{1}^{\xi _{1}}\right) \right)
^{\omega _{1}},\left( \mathbf{I}_{\otimes }\left( m,f_{2}^{\xi _{2}}\right)
\right) ^{\omega _{2}},...,\left( \mathbf{I}_{\otimes }\left( m,f_{n}^{\xi
_{n}}\right) \right) ^{\omega _{n}}\right]  \label{fff}
\end{equation}%
for all $\omega _{j},\xi _{j}\in \left( 0,\infty \right) ,$ $j=0,1,2,...n$.
\end{thrm}

\noindent \textbf{Proof.}\ Let $e\in (0,\infty ]$ be the neutral element of $%
\otimes $ and $\mathbf{I}_{\otimes }\left( m,f_{i}^{\xi _{i}}\right) =p_{i}^{%
\frac{1}{\omega _{i}}}<\infty $ for all $i=1,2,...,n$. So, for any $%
\varepsilon >0$, there exist $p_{i(\varepsilon )}^{\frac{1}{\omega _{i}}}$
such that $m\left( \left\{ f_{i}^{\xi _{i}}\geq p_{i(\varepsilon )}^{\frac{1%
}{\omega _{i}}}\right\} \right) =m\left( \left\{ f_{i}\geq p_{i(\varepsilon
)}^{\frac{1}{\xi _{i}\omega _{i}}}\right\} \right) =M\sb{i},$ where $%
p_{i(\varepsilon )}^{\frac{1}{\omega _{i}}}\otimes M\sb{i}\geq \left(
p_{i}-\varepsilon \right) ^{\frac{1}{\omega _{i}}}$ for all $i=1,2,...,n.$
The comonotonicity of $f_{1},f_{2},...,f_{n}$ and the monotonicity of $H$
imply that 
\begin{eqnarray*}
&&m\left( \left\{ H\left( f_{1},f_{2},...,f_{n}\right) \geq
H(p_{1(\varepsilon )}^{\frac{1}{\xi _{1}\omega _{1}}},p_{2(\varepsilon )}^{%
\frac{1}{\xi _{2}\omega _{2}}},...,p_{n(\varepsilon )}^{\frac{1}{\xi
_{n}\omega _{n}}})\right\} \right) \\
&\geq &m\left( \left\{ f_{1}\left( x\right) \geq p_{1(\varepsilon )}^{\frac{1%
}{\xi _{1}\omega _{1}}}\right\} \right) \wedge m\left( \left\{ f_{2}\left(
x\right) \geq p_{2(\varepsilon )}^{\frac{1}{\xi _{2}\omega _{2}}}\right\}
\right) \wedge ...\wedge m\left( \left\{ f_{n}\left( x\right) \geq
p_{n(\varepsilon )}^{\frac{1}{\xi _{n}\omega _{n}}}\right\} \right) \\
&=&M_{1}\wedge M_{2}\wedge ...\wedge M_{n}.
\end{eqnarray*}%
Since $p_{i(\varepsilon )}^{\frac{1}{\xi _{i}\omega _{i}}}\geq
p_{i(\varepsilon )}$, then we have 
\begin{eqnarray*}
&&\left[ \sup \left( t\otimes m(\{\left( H\left(
f_{1},f_{2},...,f_{n}\right) \right) ^{\xi _{0}}\geq t\})\;|\;t\in (0,\infty
])\right) \right] ^{^{\omega _{0}}} \\
&&\overset{}{\geq }\left[ \left( H(p_{1(\varepsilon )}^{\frac{1}{\xi
_{1}\omega _{1}}},p_{2(\varepsilon )}^{\frac{1}{\xi _{2}\omega _{2}}%
},...,p_{n(\varepsilon )}^{\frac{1}{\xi _{n}\omega _{n}}})\right) ^{\xi
_{0}}\otimes \left( M_{1}\wedge M_{2}\wedge ...\wedge M_{n}\right) \right]
^{^{\omega _{0}}} \\
&&\overset{}{=}\left( 
\begin{array}{c}
\left[ \left( H(p_{1(\varepsilon )}^{\frac{1}{\xi _{1}\omega _{1}}%
},p_{2(\varepsilon )}^{\frac{1}{\xi _{2}\omega _{2}}},...,p_{n(\varepsilon
)}^{\frac{1}{\xi _{n}\omega _{n}}})\right) ^{\xi _{0}}\otimes M_{1}\right]
^{\omega _{0}} \\ 
\wedge \left[ \left( H(p_{1(\varepsilon )}^{\frac{1}{\xi _{1}\omega _{1}}%
},p_{2(\varepsilon )}^{\frac{1}{\xi _{2}\omega _{2}}},...,p_{n(\varepsilon
)}^{\frac{1}{\xi _{n}\omega _{n}}})\right) ^{\xi _{0}}\otimes M_{2}\right]
^{\omega _{0}} \\ 
\wedge ...\wedge \left[ \left( H(p_{1(\varepsilon )}^{\frac{1}{\xi
_{1}\omega _{1}}},p_{2(\varepsilon )}^{\frac{1}{\xi _{2}\omega _{2}}%
},...,p_{n(\varepsilon )}^{\frac{1}{\xi _{n}\omega _{n}}})\right) ^{\xi
_{0}}\otimes M_{n}\right] ^{\omega _{0}}%
\end{array}%
\right) \\
&&\overset{}{\geq }\left( 
\begin{array}{c}
H\left( \left( p_{1(\varepsilon )}^{\frac{1}{\omega _{1}}}\otimes
M_{1}\right) ^{\omega _{1}},p_{2(\varepsilon )}^{\frac{1}{\xi _{2}\omega _{2}%
}},...,p_{n(\varepsilon )}^{\frac{1}{\xi _{n}\omega _{n}}}\right) \\ 
\wedge H\left( p_{1(\varepsilon )}^{\frac{1}{\xi _{1}\omega _{1}}},\left(
p_{2(\varepsilon )}^{\frac{1}{\omega _{2}}}\otimes M_{2}\right) ^{\omega
_{2}},p_{3(\varepsilon )}^{\frac{1}{\xi _{3}\omega _{3}}},...,p_{n(%
\varepsilon )}^{\frac{1}{\xi _{n}\omega _{n}}}\right) \\ 
\wedge ...\wedge H\left( p_{1(\varepsilon )}^{\frac{1}{\xi _{1}\omega _{1}}%
},p_{2(\varepsilon )}^{\frac{1}{\xi _{2}\omega _{2}}},...,p_{n-1(\varepsilon
)}^{\frac{1}{\xi _{n-1}\omega _{n-1}}},\left( p_{n(\varepsilon )}^{\frac{1}{%
\omega _{n}}}\otimes M_{n}\right) ^{\omega _{n}}\right)%
\end{array}%
\right) \\
&&\overset{}{\geq }\left( 
\begin{array}{c}
H\left( \left( p_{1}-\varepsilon \right) ,p_{2(\varepsilon )}^{\frac{1}{\xi
_{2}\omega _{2}}},...,p_{n(\varepsilon )}^{\frac{1}{\xi _{n}\omega _{n}}%
}\right) \\ 
\wedge H\left( p_{1(\varepsilon )}^{\frac{1}{\xi _{1}\omega _{1}}},\left(
p_{2}-\varepsilon \right) ,p_{3(\varepsilon )}^{\frac{1}{\xi _{3}\omega _{3}}%
},...,p_{n(\varepsilon )}^{\frac{1}{\xi _{n}\omega _{n}}}\right) \\ 
\wedge ...\wedge H\left( p_{1(\varepsilon )}^{\frac{1}{\xi _{1}\omega _{1}}%
},p_{2(\varepsilon )}^{\frac{1}{\xi _{2}\omega _{2}}},...,p_{n-1(\varepsilon
)}^{\frac{1}{\xi _{n-1}\omega _{n-1}}},\left( p_{n}-\varepsilon \right)
\right)%
\end{array}%
\right) \\
&&\overset{}{\geq }H\left[ (p_{1}-\varepsilon ),(p_{2}-\varepsilon
),...,(p_{n}-\varepsilon )\right] ,
\end{eqnarray*}%
whence $\left[ \mathbf{I}_{\otimes }\left( m,\left( H\left(
f_{1},f_{2},...,f_{n}\right) \right) ^{\xi _{0}}\right) \right] ^{^{\omega
_{0}}}\geq H\left[ p_{1},p_{2},...,p_{n}\right] $ follows from the
continuity of $H$ and the arbitrariness of $\varepsilon $. And the theorem
is proved. $\Box $\newline

\begin{remk}
\label{remark1} (i) If $m$ $\in \mathcal{M}_{e}^{(X,\mathcal{A})}$, then the
condition $a\otimes m\left( X\right) \leq a$ for all $a\in \lbrack 0,\infty
] $ holds readily.\newline
(ii) We can replace the condition \textquotedblleft $a\otimes m\left(
X\right) \leq a$ for all $a\in \lbrack 0,\infty ]"$ with \textquotedblleft $%
a\otimes b\leq a$ for all $a,b\in \lbrack 0,\infty ]".$
\end{remk}

\begin{corl}
\label{th3-1} Let $f,g\in $ $\mathcal{F}^{(X,\mathcal{A})}$ be two
comonotone measurable functions and $\otimes \colon \lbrack 0,\infty ]\sp{2}%
\rightarrow \lbrack 0,\infty ]$ be a smallest pseudo-multiplication on $%
[0,\infty ]$ with neutral element $e\in (0,\infty ]$ and $m$ $\in \mathcal{M}%
^{(X,\mathcal{A})}$ be a monotone measure such that $a\otimes m\left(
X\right) \leq a$ for all $a\in \lbrack 0,\infty ],$ $\mathbf{I}_{\otimes
}\left( m,g^{\xi _{2}}\right) <\infty $ and $\mathbf{I}_{\otimes }\left(
m,f^{\xi _{1}}\right) $ $<\infty $. Let $\star \colon \lbrack 0,\infty )%
\sp{2}\rightarrow \lbrack 0,\infty )$ be continuous and nondecreasing in
both arguments. If 
\begin{equation}
\left[ \left( t_{1}\star t_{2}\right) ^{\xi _{0}}\otimes c\right] ^{\omega
_{0}}\geq \left[ \left( t_{1}^{\xi _{1}}\otimes c\right) ^{\omega _{1}}\star
t_{2}\right] \vee \left[ t_{1}\star (t_{2}^{\xi _{2}}\otimes c)^{\omega _{2}}%
\right] ,  \label{poio}
\end{equation}%
then the inequality 
\begin{equation}
\left[ \mathbf{I}_{\otimes }\left( m,(f\star g)^{\xi _{0}}\right) \right]
^{\omega _{0}}\geq \left[ \mathbf{I}_{\otimes }\left( m,f^{\xi _{1}}\right) %
\right] ^{\omega _{1}}\star \left[ \mathbf{I}_{\otimes }\left( m,g^{\xi
_{2}}\right) \right] ^{\omega _{2}}  \label{eq3-2}
\end{equation}%
holds, where $x^{\frac{1}{\xi _{i}\omega _{i}}}\geq x$ for all $x\in \lbrack
0,\infty ),i=1,2$ and $\omega _{j},\xi _{j}\in \left( 0,\infty \right)
,j=0,1,2$.
\end{corl}

The following example shows that the condition of $x^{\frac{1}{\xi
_{i}\omega _{i}}}\geq x$ for all $x\in \lbrack 0,\infty )$ and $i=1,2$ in
Corollary \ref{th3-1} (and thus the condition $x^{\frac{1}{\xi _{i}\omega
_{i}}}\geq x$ for all $x\in \lbrack 0,\infty )$ and $i=1,2,...n$ in Theorem %
\ref{th3-2}) is inevitable.

\begin{exa}
Let $X=[0,1],\star =\wedge ,\xi _{0}=\omega _{0}=1,\xi _{i}=\frac{1}{2}%
,\omega _{i}=1$ for $i=1,2$. Let $f\left( x\right) =x,g\left( x\right) $ $=1$
for all $x\in \lbrack 0,1]$ and the monotone measure $m$ be the Lebesgue
measure. If $\ \otimes \colon \lbrack 0,1]\sp{2}\rightarrow \lbrack 0,1]$ is
minimum (i.e., for Sugeno integral), then (\ref{poio}) holds readily for all 
$t_{1},t_{2},c\in \lbrack 0,1]$ and a straightforward calculus shows that 
\begin{eqnarray*}
&&\underset{}{\left( i\right) }\ \mathbf{I}_{Min}\left( m,f^{\frac{1}{2}%
}\right) \underset{}{=}\mathbf{Su}\left( m,f^{\frac{1}{2}}\right) \underset{}%
{=}\bigvee_{\alpha \in \lbrack 0,\allowbreak 1]}\left[ \alpha \wedge m\left(
\left\{ \sqrt{x}\geq \alpha \right\} \right) \right] \underset{}{=}\frac{1}{2%
}\left( \sqrt{5}-1\right) , \\
&&\underset{}{\left( ii\right) \text{ }}\mathbf{I}_{Min}\left( m,g^{\frac{1}{%
2}}\right) \underset{}{=}\mathbf{Su}\left( m,g^{\frac{1}{2}}\right) \underset%
{}{=}1, \\
&&\underset{}{\left( iii\right) }\text{ }\mathbf{I}_{Min}\left( m,(f\wedge
g)\right) \underset{}{=}\mathbf{Su}\left( m,f\right) \underset{}{=}%
\bigvee_{\alpha \in \lbrack 0,\allowbreak 1]}\left[ \alpha \wedge m\left(
\left\{ x\geq \alpha \right\} \right) \right] \underset{}{=}\frac{1}{2}.
\end{eqnarray*}%
Therefore: 
\begin{eqnarray*}
\left[ \mathbf{I}_{\otimes }\left( m,(f\star g)^{\xi _{0}}\right) \right]
^{\omega _{0}} &=&\text{ }\mathbf{I}_{Min}\left( m,(f\wedge g)\right) =\frac{%
1}{2}<\left[ \mathbf{I}_{\otimes }\left( m,f^{\xi _{1}}\right) \right]
^{\omega _{1}}\star \left[ \mathbf{I}_{\otimes }\left( m,g^{\xi _{2}}\right) %
\right] ^{\omega _{2}} \\
&=&\mathbf{I}_{Min}\left( m,f^{\frac{1}{2}}\right) \wedge \mathbf{I}%
_{Min}\left( m,g^{\frac{1}{2}}\right) =\frac{1}{2}\left( \sqrt{5}-1\right) ,
\end{eqnarray*}%
which violates Corollary \ref{th3-1}.
\end{exa}

\begin{remk}
\label{remark3.5} If $(x\star e)\vee (e\star x)\leq x$ and $x^{\frac{1}{\xi
_{0}\omega _{0}}}\leq x\leq x^{\frac{1}{\xi _{i}\omega _{i}}},$ $x\geq x^{%
\frac{\omega _{i}}{\omega _{_{0}}}}$for any $x\in \lbrack 0,\infty )$ and $%
\omega _{i},\xi _{i}\in \left( 0,\infty \right) ,i=1,2$ and $\left( .\right)
^{\omega _{0}}\;$is superdistributive over $\otimes $ and $\left( .\right)
^{\omega _{i}},i=1,2$ are subdistributive over $\otimes $ and $\otimes $
dominates $\star $, then (\ref{poio}) holds readily. Indeed,%
\begin{eqnarray*}
&&\left[ \left( t_{1}\star t_{2}\right) ^{\xi _{0}}\otimes c\right] ^{\omega
_{0}}\overset{}{\geq }\left( t_{1}\star t_{2}\right) ^{\omega _{0}\xi
_{0}}\otimes c^{\omega _{0}}\overset{}{\geq }\left[ \left( t_{1}\star
t_{2}\right) \otimes c^{\omega _{1}}\right] \\
&&\overset{}{\geq }\left[ \left( t_{1}\star t_{2}\right) \otimes \left(
c^{\omega _{1}}\star e\right) \right] \overset{}{\geq }\left[ \left(
t_{1}\otimes c^{\omega _{1}}\right) \star \left( t_{2}\otimes e\right) %
\right] \\
&&\overset{}{=}\left[ \left( t_{1}\otimes c^{\omega _{1}}\right) \star t_{2}%
\right] \overset{}{\geq }\left[ \left( t_{1}^{\omega _{1}\xi _{1}}\otimes
c^{\omega _{1}}\right) \star t_{2}\right] \geq \left[ \left( t_{1}^{\xi
_{1}}\otimes c\right) ^{\omega _{1}}\star t_{2}\right] ,
\end{eqnarray*}%
and $\left[ \left( t_{1}\star t_{2}\right) ^{\xi _{0}}\otimes c\right]
^{\omega _{0}}\geq \left[ t_{1}\star (t_{2}^{\xi _{2}}\otimes c)^{\omega
_{2}}\right] $ follows similarly, i.e., 
\begin{eqnarray*}
&&\left[ \left( t_{1}\star t_{2}\right) ^{\xi _{0}}\otimes c\right] ^{\omega
_{0}}\overset{}{\geq }\left( t_{1}\star t_{2}\right) ^{\omega _{0}\xi
_{0}}\otimes c^{\omega _{0}}\overset{}{\geq }\left[ \left( t_{1}\star
t_{2}\right) \otimes c^{\omega _{1}}\right] \\
&&\overset{}{\geq }\left[ \left( t_{1}\star t_{2}\right) \otimes \left(
c^{\omega _{2}}\star e\right) \right] \overset{}{\geq }\left[ \left(
t_{1}\otimes e\right) \star \left( t_{2}\otimes c^{\omega _{2}}\right) %
\right] \\
&&\overset{}{=}\left[ t_{1}\star \left( t_{2}\otimes c^{\omega _{2}}\right) %
\right] \overset{}{\geq }\left[ t_{1}\star \left( t_{2}^{\omega _{2}\xi
_{2}}\otimes c^{\omega _{2}}\right) \right] \geq t_{1}\star \left(
t_{2}^{\xi _{2}}\otimes c\right) ^{\omega _{2}}.
\end{eqnarray*}
\end{remk}

We get an inequality related to the H\"{o}lder type inequality whenever $\xi
_{0}=\omega _{0}=1,\xi _{1}=p,\omega _{1}=\frac{1}{p},\xi _{2}=q$ and $%
\omega _{2}=\frac{1}{q}$ for all $p,q\in \left( 0,\infty \right) .$

\begin{corl}
Let $f,g\in $ $\mathcal{F}^{(X,\mathcal{A})}$ be two comonotone measurable
functions and $\otimes \colon \lbrack 0,\infty ]\sp{2}\rightarrow \lbrack
0,\infty ]$ be a smallest pseudo-multiplication on $[0,\infty ]$ with
neutral element $e\in (0,\infty ]$ and $m$ $\in \mathcal{M}^{(X,\mathcal{A}%
)} $ be a monotone measure such that $a\otimes m\left( X\right) \leq a$ for
all $a\in \lbrack 0,\infty ],$ $\mathbf{I}_{\otimes }\left( m,g^{q}\right)
<\infty $ and $\mathbf{I}_{\otimes }\left( m,f^{p}\right) $ $<\infty $. Let $%
\star \colon \lbrack 0,\infty )\sp{2}\rightarrow \lbrack 0,\infty )$ be
continuous and nondecreasing in both arguments. If 
\begin{equation*}
\left[ \left( a\star b\right) \otimes c\right] \geq \left[ \left(
a^{p}\otimes c\right) ^{\frac{1}{p}}\star b\right] \vee \left[ a\star
(b^{q}\otimes c)^{\frac{1}{q}}\right] ,
\end{equation*}%
then the inequality 
\begin{equation*}
\left[ \mathbf{I}_{\otimes }\left( m,(f\star g)\right) \right] \geq \left[ 
\mathbf{I}_{\otimes }\left( m,f^{p}\right) \right] ^{\frac{1}{p}}\star \left[
\mathbf{I}_{\otimes }\left( m,g^{q}\right) \right] ^{\frac{1}{q}}
\end{equation*}%
holds for all $p,q\in \left( 0,\infty \right) $.
\end{corl}

Again, we get an inequality related to the Minkowski type whenever $\xi
_{0}=\xi _{1}=\xi _{2}=s$ and $\omega _{0}=\omega _{1}=\omega _{2}=\frac{1}{s%
}$ for all $s\in \left( 0,\infty \right) .$

\begin{corl}
Let $f,g\in $ $\mathcal{F}^{(X,\mathcal{A})}$ be two comonotone measurable
functions and $\otimes \colon \lbrack 0,\infty ]\sp{2}\rightarrow \lbrack
0,\infty ]$ be a smallest pseudo-multiplication on $[0,\infty ]$ with
neutral element $e\in (0,\infty ]$ and $m$ $\in \mathcal{M}^{(X,\mathcal{A}%
)} $ be a monotone measure such that $a\otimes m\left( X\right) \leq a$ for
all $a\in \lbrack 0,\infty ],$ $\mathbf{I}_{\otimes }\left( m,f^{s}\right)
<\infty $ and $\mathbf{I}_{\otimes }\left( m,g^{s}\right) <\infty .$ Let $%
\star \colon \lbrack 0,\infty )\sp{2}\rightarrow \lbrack 0,\infty )$ be
continuous and nondecreasing in both arguments. If 
\begin{equation}
\left[ \left( a\star b\right) ^{s}\otimes c)\right] ^{\frac{1}{s}}\geq \left[
\left( a^{s}\otimes c\right) ^{\frac{1}{s}}\star b\right] \vee \left[ a\star
(b^{s}\otimes c)^{\frac{1}{s}}\right] ,  \label{jjj}
\end{equation}%
then the inequality 
\begin{equation*}
\left( \mathbf{I}_{\otimes }\left( m,\left( f\star g\right) ^{s}\right)
\right) ^{\frac{1}{s}}\geq \left( \mathbf{I}_{\otimes }\left( m,f^{s}\right)
\right) ^{\frac{1}{s}}\star \left( \mathbf{I}_{\otimes }\left(
m,g^{s}\right) \right) ^{\frac{1}{s}}
\end{equation*}%
holds for all $s>0.$
\end{corl}

Specially, when $s=1$ we have the Chebyshev inequality.

\begin{corl}
Let $f,g\in $ $\mathcal{F}^{(X,\mathcal{A})}$ be two comonotone measurable
functions and $\otimes \colon \lbrack 0,\infty ]\sp{2}\rightarrow \lbrack
0,\infty ]$ be a smallest pseudo-multiplication on $[0,\infty ]$ with
neutral element $e\in (0,\infty ]$ and $m$ $\in \mathcal{M}^{(X,\mathcal{A}%
)} $ be a monotone measure such that $a\otimes m\left( X\right) \leq a$ for
all $a\in \lbrack 0,\infty ],$ $\mathbf{I}_{\otimes }\left( m,f\right)
<\infty $ and $\mathbf{I}_{\otimes }\left( m,g\right) <\infty $. Let $\star
\colon \lbrack 0,\infty )\sp{2}\rightarrow \lbrack 0,\infty )$ be continuous
and nondecreasing in both arguments. If 
\begin{equation*}
\left( a\star b\right) \otimes c)\geq \left[ \left( a\otimes c\right) \star b%
\right] \vee \left[ a\star (b\otimes c)\right] ,
\end{equation*}%
then the inequality 
\begin{equation*}
\mathbf{I}_{\otimes }\left( m,\left( f\star g\right) \right) \geq \mathbf{I}%
_{\otimes }\left( m,f\right) \star \mathbf{I}_{\otimes }\left( m,g\right)
\end{equation*}%
holds.
\end{corl}

\begin{remk}
If $\otimes $ is minimum (i.e., for Sugeno integral) and $n$-place function $%
H:[0,\infty )^{n}\rightarrow \lbrack 0,\infty )$ is continuous and
nondecreasing and bounded from above by minimum, then (\ref{jjj}) holds
readily whenever $x^{\frac{1}{\xi _{0}\omega _{0}}}\leq x\leq x^{\frac{1}{%
\xi _{i}\omega _{i}}}$ and $x\geq x^{\frac{\omega _{i}}{\omega _{0}}}$ for
all $x\in \lbrack 0,\infty )$ and $\omega _{i},\xi _{i}\in \left( 0,\infty
\right) ,i=1,2,...n$. Indeed,%
\begin{eqnarray*}
&&\left[ H^{\xi _{0}}\left( p_{1},p_{2},...,p_{n}\right) \wedge c\right]
^{\omega _{0}}=\left[ H^{\omega _{0}\xi _{0}}\left(
p_{1},p_{2},...,p_{n}\right) \wedge c^{\omega _{0}}\right] \overset{}{\geq }%
\left[ H\left( p_{1},p_{2},...,p_{n}\right) \wedge c^{\omega _{0}}\right] \\
&\geq &\left[ H\left( p_{1}^{\omega _{1}\xi _{1}},p_{2},...,p_{n}\right)
\wedge c^{\omega _{1}}\right] \geq \left[ H\left( p_{1}^{\omega _{1}\xi
_{1}},p_{2},...,p_{n}\right) \wedge (p_{1}^{\omega _{1}\xi _{1}}\wedge
c^{\omega _{1}})\right] \\
&\geq &\left[ H\left( p_{1}^{\omega _{1}\xi _{1}},p_{2},...,p_{n}\right)
\wedge \left( (p_{1}^{\omega _{1}\xi _{1}}\wedge c^{\omega _{1}})\wedge
p_{2}\wedge ...\wedge p_{n}\right) \right] \\
&\geq &\left[ H\left( p_{1}^{\omega _{1}\xi _{1}},p_{2},...,p_{n}\right)
\wedge H\left( (p_{1}^{\omega _{1}\xi _{1}}\wedge c^{\omega
_{1}}),p_{2},...,p_{n}\right) \right] \\
&\geq &\left[ H\left( (p_{1}^{\omega _{1}\xi _{1}}\wedge c^{\omega
_{1}}),p_{2},...,p_{n}\right) \wedge H\left( (p_{1}^{\omega _{1}\xi
_{1}}\wedge c^{\omega _{1}}),p_{2},...,p_{n}\right) \right] \\
&=&H\left( (p_{1}^{\omega _{1}\xi _{1}}\wedge c^{\omega
_{1}}),p_{2},...,p_{n}\right) .
\end{eqnarray*}%
and the others follow similarly. Thus the following results hold.
\end{remk}

\begin{corl}
Let $n$-place function $H:[0,\infty )^{n}\rightarrow \lbrack 0,\infty )$ be
continuous and nondecreasing and bounded from above by minimum. Then for any
comontone system $f_{1},f_{2},...,f_{n}\in $ $\mathcal{F}^{(X,\mathcal{A})}$
and a monotone measure $m$ $\in \mathcal{M}^{(X,\mathcal{A})}$ such that $%
\mathbf{Su}\left( m,f_{i}^{\xi _{i}}\right) <\infty ,x^{\frac{1}{\xi
_{0}\omega _{0}}}\leq x\leq x^{\frac{1}{\xi _{i}\omega _{i}}}$ and $x\geq x^{%
\frac{\omega _{i}}{\omega _{0}}}$ for all $x\in \lbrack 0,\infty )$ and $%
\omega _{i},\xi _{i}\in \left( 0,\infty \right) ,i=1,2,...n$, it holds%
\begin{equation*}
\left[ \mathbf{Su}\left( m,\left( H\left( f_{1},...,f_{n}\right) \right)
^{\xi _{0}}\right) \right] ^{^{\omega _{0}}}\overset{}{\geq }H\left[ \left( 
\mathbf{Su}\left( m,f_{1}^{\xi _{1}}\right) \right) ^{\omega _{1}},\left( 
\mathbf{Su}\left( m,f_{2}^{\xi _{2}}\right) \right) ^{\omega
_{2}},...,\left( \mathbf{Su}\left( m,f_{n}^{\xi _{n}}\right) \right)
^{\omega _{n}}\right]
\end{equation*}%
for all $\omega _{j},\xi _{j}\in \left( 0,\infty \right) ,$ $j=0,1,2,...n$.
\end{corl}

\begin{corl}
Let $f_{1},f_{2}\in $ $\mathcal{F}^{(X,\mathcal{A})}$ be two comonotone
measurable functions. Let $\star \colon \lbrack 0,\infty )\sp{2}\rightarrow
\lbrack 0,\infty )$ be continuous and nondecreasing in both arguments and
bounded from above by minimum and $m\in \mathcal{M}^{(X,\mathcal{A})}$ be a
monotone measure such that $\mathbf{Su}\left( m,f_{i}^{\xi _{i}}\right)
<\infty ,x^{\frac{1}{\xi _{0}\omega _{0}}}\leq x\leq x^{\frac{1}{\xi
_{i}\omega _{i}}}$ and $x\geq x^{\frac{\omega _{i}}{\omega _{0}}}$ for all $%
x\in \lbrack 0,\infty )$ and $\omega _{i},\xi _{i}\in \left( 0,\infty
\right) ,i=1,2$, it holds%
\begin{equation*}
\left[ \mathbf{Su}\left( m,\left( f_{1}\star f_{2}\right) ^{\xi _{0}}\right) %
\right] ^{^{\omega _{0}}}\geq \left[ \mathbf{Su}\left( m,f_{1}^{\xi
_{1}}\right) \right] ^{\omega _{1}}\star \left[ \mathbf{Su}\left(
m,f_{2}^{\xi _{2}}\right) \right] ^{\omega _{2}}
\end{equation*}%
for all $\omega _{j},\xi _{j}\in \left( 0,\infty \right) ,$ $j=0,1,2$.
\end{corl}

\begin{corl}
(\textrm{\cite{OuyMesAg}) }Let $f,g\in $ $\mathcal{F}^{(X,\mathcal{A})} $ be
two comonotone measurable functions. Let $\star \colon \lbrack 0,\infty )%
\sp{2}\rightarrow \lbrack 0,\infty )$ be continuous and nondecreasing in
both arguments and bounded from above by minimum and $m\in \mathcal{M}^{(X,%
\mathcal{A})}$ be a monotone measure such that $\mathbf{Su}\left(
m,f^{s}\right) <\infty ,\mathbf{Su}\left( m,g^{s}\right) <\infty $. Then the
inequality 
\begin{equation*}
\left[ \mathbf{Su}\left( m,\left( f\star g\right) ^{s}\right) \right] ^{%
\frac{1}{s}}\geq \left[ \mathbf{Su}\left( m,f^{s}\right) \right] ^{\frac{1}{s%
}}\star \left[ \mathbf{Su}\left( m,g^{s}\right) \right] ^{\frac{1}{s}}
\end{equation*}%
holds for all $0<s<\infty $.
\end{corl}

\begin{corl}
Let $f,g\in $ $\mathcal{F}^{(X,\mathcal{A})}$ be two comonotone measurable
functions. Let $\star \colon \lbrack 0,\infty )\sp{2}\rightarrow \lbrack
0,\infty )$ be continuous and nondecreasing in both arguments and bounded
from above by minimum and $m\in \mathcal{M}^{(X,\mathcal{A})}$ be a monotone
measure such that $\mathbf{Su}\left( m,f^{p}\right) <\infty ,\mathbf{Su}%
\left( m,g^{q}\right) <\infty $. Then the inequality 
\begin{equation*}
\mathbf{Su}\left( m,\left( f\star g\right) \right) \geq \left[ \mathbf{Su}%
\left( m,f^{p}\right) \right] ^{\frac{1}{p}}\star \left[ \mathbf{Su}\left(
m,g^{q}\right) \right] ^{\frac{1}{q}}
\end{equation*}%
holds, where $x\geq x^{\frac{1}{p}},x\geq x^{\frac{1}{q}}$for all $x\in
\lbrack 0,\infty )$ and $p,q\in (0,\infty )$.
\end{corl}

\begin{corl}
(\cite{MesOuy09}) Let $f,g\in $ $\mathcal{F}^{(X,\mathcal{A})}$ be two
comonotone measurable functions. Let $\star \colon \lbrack 0,\infty )\sp{2}%
\rightarrow \lbrack 0,\infty )$ be continuous and nondecreasing in both
arguments and bounded from above by minimum and $m\in \mathcal{M}^{(X,%
\mathcal{A})}$ be a monotone measure such that $\mathbf{Su}\left( m,f\right)
<\infty ,\mathbf{Su}\left( m,g\right) <\infty $. Then the inequality 
\begin{equation*}
\mathbf{Su}\left( m,f\star g\right) \geq \mathbf{Su}\left( m,f\right) \star 
\mathbf{Su}\left( m,g\right)
\end{equation*}%
holds.
\end{corl}

Notice that when working on $[0,1]$ in Theorem \ref{th3-1}, we mostly deal
with $e=1$, then $\otimes =\circledast $ is semicopula (t-seminorm) and the
following results hold.

\begin{corl}
Let a non-decreasing $n$-place function $H:[0,\infty )^{n}\rightarrow
\lbrack 0,\infty )$ such that $H$ be continuous. If semicopula $\circledast $
satisfies 
\begin{eqnarray*}
&&\left[ \left( H\left( p_{1},p_{2},...,p_{n}\right) \right) ^{\xi
_{0}}\circledast c\right] ^{\omega _{0}}\overset{}{\geq }H\left( \left(
p_{1}^{\xi _{1}}\circledast c\right) ^{\omega _{1}},p_{2},...,p_{n}\right)
\vee \\
&&H\left( p_{1},\left( p_{2}^{\xi _{2}}\circledast c\right) ^{\omega
_{2}},p_{3},...,p_{n}\right) \vee ...\vee H\left(
p_{1},p_{2},...,p_{n-1},\left( p_{n}^{\xi _{n}}\circledast c\right) ^{\omega
_{n}}\right) ,
\end{eqnarray*}%
then for any comontone system $f_{1},f_{2},...,f_{n}\in $ $\mathcal{F}%
_{1}^{(X,\mathcal{A})}$ and a monotone measure $m$ $\in \mathcal{M}_{1}^{(X,%
\mathcal{A})}$, it holds%
\begin{equation*}
\left[ \mathbf{I}_{\circledast }\left( m,\left( H\left(
f_{1},...,f_{n}\right) \right) ^{\xi _{0}}\right) \right] ^{^{\omega _{0}}}%
\overset{}{\geq }H\left[ \left( \mathbf{I}_{\circledast }\left( m,f_{1}^{\xi
_{1}}\right) \right) ^{\omega _{1}},\left( \mathbf{I}_{\circledast }\left(
m,f_{2}^{\xi _{2}}\right) \right) ^{\omega _{2}},...,\left( \mathbf{I}%
_{\circledast }\left( m,f_{n}^{\xi _{n}}\right) \right) ^{\omega _{n}}\right]
,
\end{equation*}%
where $\omega _{i}\xi _{i}\geq 1$ for all $\omega _{j},\xi _{j}\in \left(
0,\infty \right) ,i=1,2,...n$ and $j=0,1,2,...n$.
\end{corl}

\begin{corl}
\label{colr3-1} Let $f,g\in $ $\mathcal{F}_{\left[ 0,1\right] }^{(X,\mathcal{%
A})}$ be two comonotone measurable functions. Let $\star \colon \lbrack 0,1]%
\sp{2}\rightarrow \lbrack 0,1]$ be continuous and nondecreasing in both
arguments. If semicopula $\circledast $ satisfies 
\begin{equation}
\left[ \left( a\star b\right) ^{\alpha }\circledast c\right] ^{\lambda }\geq %
\left[ \left( a^{\beta }\circledast c\right) ^{\upsilon }\star b\right] \vee %
\left[ a\star (b^{\gamma }\circledast c)^{\tau }\right] ,  \label{mmm}
\end{equation}%
then the inequality 
\begin{equation*}
\left[ \mathbf{I}_{\circledast }\left( m,(f\star g)^{\alpha }\right) \right]
^{\lambda }\geq \left[ \mathbf{I}_{\circledast }\left( m,f^{\beta }\right) %
\right] ^{\upsilon }\star \left[ \mathbf{I}_{\circledast }\left( m,g^{\gamma
}\right) \right] ^{\tau }
\end{equation*}%
holds for all $\alpha ,\beta ,\gamma ,\lambda ,\upsilon ,\tau \in \left(
0,\infty \right) ,\gamma \tau \geq 1,\beta \upsilon \geq 1$ and for any $%
m\in \mathcal{M}_{1}^{(X,\mathcal{A})}$.
\end{corl}

Let $\alpha =\beta =\gamma =s$ and $\lambda =\upsilon =\tau =\frac{1}{s}$
for all $s\in \left( 0,\infty \right) $, then we get the reverse Minkowski
type inequality for seminormed fuzzy integrals.

\begin{corl}
Let $f,g\in $ $\mathcal{F}_{\left[ 0,1\right] }^{(X,\mathcal{A})}$ be two
comonotone measurable functions. Let $\star \colon \lbrack 0,1]\sp{2}%
\rightarrow \lbrack 0,1]$ be continuous and nondecreasing in both arguments.
If semicopula $\circledast $ satisfies 
\begin{equation*}
\left[ \left( a\star b\right) ^{s}\circledast c)\right] ^{\frac{1}{s}}\geq %
\left[ \left( a^{s}\circledast c\right) ^{\frac{1}{s}}\star b\right] \vee %
\left[ a\star (b^{s}\circledast c)^{\frac{1}{s}}\right] ,
\end{equation*}%
then the inequality 
\begin{equation*}
\left( \mathbf{I}_{\circledast }\left( m,\left( f\star g\right) ^{s}\right)
\right) ^{\frac{1}{s}}\geq \left( \mathbf{I}_{\circledast }\left(
m,f^{s}\right) \right) ^{\frac{1}{s}}\star \left( \mathbf{I}_{\circledast
}\left( m,g^{s}\right) \right) ^{\frac{1}{s}}
\end{equation*}%
holds for any $m\in \mathcal{M}_{1}^{(X,\mathcal{A})}$ and for all $%
0<s<\infty $.
\end{corl}

Again, we get the Chebyshev type inequality for seminormed fuzzy integrals
whenever $s=1$ \cite{OuyMes09}.

\begin{corl}
\label{col3-1}Let $f,g\in $ $\mathcal{F}_{\left[ 0,1\right] }^{(X,\mathcal{A}%
)}$ be two comonotone measurable functions. Let $\star \colon \lbrack 0,1]%
\sp{2}\rightarrow \lbrack 0,1]$ be continuous and nondecreasing in both
arguments. If semicopula $\circledast $ satisfies 
\begin{equation*}
\left[ \left( a\star b\right) \circledast c)\right] \geq \left[ \left(
a\circledast c\right) \star b\right] \vee \left[ a\star (b\circledast c)%
\right] ,
\end{equation*}%
then the inequality 
\begin{equation*}
\mathbf{I}_{\circledast }\left( m,\left( f\star g\right) \right) \geq 
\mathbf{I}_{\circledast }\left( m,f\right) \star \mathbf{I}_{\circledast
}\left( m,g\right)
\end{equation*}%
holds for any $m\in \mathcal{M}_{1}^{(X,\mathcal{A})}$.
\end{corl}

\begin{remk}
We can use an example in \cite{OuyMes09} to show that the condition of $%
\left[ \left( a\star b\right) \circledast c\right] \geq \left[ \left(
a\circledast c\right) \star b\right] \vee \left[ a\star (b\circledast c)%
\right] $ in Corollary \ref{col3-1} (and thus in Theorem \ref{th3-2}) cannot
be abandoned, and so we omit it here.
\end{remk}

Suppose the semicopula $\circledast $ further satisfies monotonicity and
associativity (i.e., it is a $t$-norm). Then, we have the following result:

\begin{corl}
Let $f,g\in $ $\mathcal{F}_{\left[ 0,1\right] }^{(X,\mathcal{A})}$ be two
comonotone measurable functions. Let $\star \colon \lbrack 0,1]\sp{2}%
\rightarrow \lbrack 0,1]$ be continuous and nondecreasing in both arguments.
If semicopula $\circledast $ be a continuous $t$-norm, then 
\begin{equation*}
\left[ \mathbf{I}_{\circledast }\left( m,(f\circledast g)^{\alpha }\right) %
\right] ^{\lambda }\geq \left( \left[ \mathbf{I}_{\circledast }\left(
m,f^{\beta }\right) \right] ^{\upsilon }\circledast \left[ \mathbf{I}%
_{\circledast }\left( m,g^{\gamma }\right) \right] ^{\tau }\right)
\end{equation*}%
holds for any $m\in \mathcal{M}_{1}^{(X,\mathcal{A})}$ and for all $\alpha
,\beta ,\gamma ,\lambda ,\upsilon ,\tau \in \left( 0,\infty \right)
,0<\alpha \lambda \leq 1,1\leq \beta \upsilon <\infty ,1\leq \gamma \tau
<\infty ,\lambda \leq \tau ,\upsilon \ $and $\alpha \leq \beta ,\gamma ,$
where $\left( .\right) ^{\alpha }$ is superdistributive over $\circledast $, 
$\circledast ^{\lambda }$ dominates $\circledast $ and $(f\circledast
g)(x)=f(x)\circledast g(x)$ for any $x\in X$.
\end{corl}

Let $\alpha =\beta =\gamma =\lambda =\upsilon =\tau =1$, then $\circledast $
is obviously dominated by itself and we have the following result:

\begin{corl}
Let $f,g\in $ $\mathcal{F}_{\left[ 0,1\right] }^{(X,\mathcal{A})}$ be two
comonotone measurable functions. Let $\star \colon \lbrack 0,1]\sp{2}%
\rightarrow \lbrack 0,1]$ be continuous and nondecreasing in both arguments.
If semicopula $\circledast $ be a continuous $t$-norm, then 
\begin{equation*}
\mathbf{I}_{\circledast }\left( m,(f\circledast g)\right) \geq \left( 
\mathbf{I}_{\circledast }\left( m,f\right) \circledast \mathbf{I}%
_{\circledast }\left( m,g\right) \right)
\end{equation*}%
holds for any $m\in \mathcal{M}_{1}^{(X,\mathcal{A})}$ and $(f\circledast
g)(x)=f(x)\circledast g(x)$ for any $x\in X$.
\end{corl}

Notice that if the semicopula ($t$-seminorm) $\circledast $ is minimum
(i.e., for Sugeno integral) and $\star $ is bounded from above by minimum,
then $\star $ is dominated by minimum. Thus the following result holds.

\begin{corl}
Let $f,g\in $ $\mathcal{F}_{\left[ 0,1\right] }^{(X,\mathcal{A})}$ be two
comonotone measurable functions. Let $\star \colon \lbrack 0,1]\sp{2}%
\rightarrow \lbrack 0,1]$ be continuous and nondecreasing in both arguments
and bounded from above by minimum. Then the inequality 
\begin{equation*}
\left[ \mathbf{Su}\left( m,\left( f\star g\right) ^{\alpha }\right) \right]
^{\lambda }\geq \left[ \mathbf{Su}\left( m,f^{\beta }\right) \right]
^{\upsilon }\star \left[ \mathbf{Su}\left( m,g^{\gamma }\right) \right]
^{\tau }
\end{equation*}%
holds for any $m\in \mathcal{M}_{1}^{(X,\mathcal{A})}$ and for all $\alpha
,\beta ,\gamma ,\lambda ,\upsilon ,\tau \in \left( 0,\infty \right)
,0<\alpha \lambda \leq 1,\beta \upsilon \geq 1,\gamma \tau \geq 1,\lambda
\leq \tau ,\upsilon .$
\end{corl}

\begin{thrm}
\label{th3-3}Let $f\in $ $\mathcal{F}^{(X,\mathcal{A})}$ be a measurable
function and $\otimes \colon \lbrack 0,\infty ]\sp{2}\rightarrow \lbrack
0,\infty ]$ be the pseudo-multiplication with neutral element $e\in
(0,\infty ]$ and $m$ $\in \mathcal{M}^{(X,\mathcal{A})}$ be a monotone
measure such that $\mathbf{I}_{\otimes }\left( m,\varphi _{2}\left( f\right)
\right) $ is finite. Let $\varphi _{i}:[0,\infty )\rightarrow \lbrack
0,\infty ),i=1,2$ be continuous strictly increasing functions. If $\ $ 
\begin{equation*}
\varphi _{1}^{-1}\left( \varphi _{1}\left( a\right) \otimes c\right) \geq
\varphi _{2}^{-1}\left( \varphi _{2}\left( a\right) \otimes c\right) ,
\end{equation*}%
then the inequality 
\begin{equation*}
\varphi _{1}^{-1}\left( \mathbf{I}_{\otimes }\left( m,\varphi _{1}\left(
f\right) \right) \right) \geq \varphi _{2}^{-1}\left( \mathbf{I}_{\otimes
}\left( m,\varphi _{2}\left( f\right) \right) \right)
\end{equation*}%
holds.
\end{thrm}

\noindent\textbf{Proof.}\ Let $e\in (0,\infty ]$ be the neutral element of $\otimes $
and $\mathbf{I}_{\otimes }\left( m,\varphi _{2}\left( f\right) \right)
=p<\infty $. So, for any $\varepsilon >0$, there exists $p_{\varepsilon }$
such that $m(\left\{ \varphi _{2}\left( f\right) \geq p_{\varepsilon
}\right\} )=M,$where $p_{\varepsilon }\otimes M\geq p-\varepsilon $. Hence, 
\begin{eqnarray*}
&&\varphi _{1}^{-1}\left( \mathbf{I}_{\otimes }\left( m,\varphi _{1}\left(
f\right) \right) \right) \overset{}{\geq }\varphi _{1}^{-1}\left( \left[
\varphi _{1}\left( \varphi _{2}^{-1}\left( p_{\varepsilon }\right) \right)
\otimes m(\{\varphi _{1}\left( f\right) \geq \varphi _{1}\left( \varphi
_{2}^{-1}\left( p_{\varepsilon }\right) \right) \})\right] \right) \\
&&\overset{}{=}\varphi _{1}^{-1}\left( \left[ \varphi _{1}\left( \varphi
_{2}^{-1}\left( p_{\varepsilon }\right) \right) \otimes m(\{\varphi
_{2}\left( f\right) \geq p_{\varepsilon }\})\right] \right) \\
&&\overset{}{\geq }\varphi _{2}^{-1}\left( \left[ \varphi _{2}\left( \varphi
_{2}^{-1}\left( p_{\varepsilon }\right) \right) \otimes m(\{\varphi
_{2}\left( f\right) \geq p_{\varepsilon }\})\right] \right) \\
&&\overset{}{=}\varphi _{2}^{-1}\left( \left[ p_{\varepsilon }\otimes M%
\right] \right) \overset{}{\geq }\varphi _{2}^{-1}\left( p-\varepsilon
\right)
\end{eqnarray*}%
whence $\varphi _{1}^{-1}\left( \mathbf{I}_{\otimes }\left( m,\varphi
_{1}\left( f\right) \right) \right) \geq \varphi _{2}^{-1}\left( p\right) $
follows from the continuity of $\varphi _{2}$ and the arbitrariness of $%
\varepsilon $. And the theorem is proved. $\Box $

If we take $\varphi _{2}\left( x\right) =x$ in Theorem \ref{th3-3}, then the
the following Jensen inequality for universal integral is recaptured.

\begin{corl}
\label{cor3.7} Let $f\in $ $\mathcal{F}^{(X,\mathcal{A})}$ be a measurable
function and $\otimes \colon \lbrack 0,\infty ]\sp{2}\rightarrow \lbrack
0,\infty ]$ be the pseudo-multiplication with neutral element $e\in
(0,\infty ]$ and $m$ $\in \mathcal{M}^{(X,\mathcal{A})}$ be a monotone
measure such that $\mathbf{I}_{\otimes }\left( m,f\right) $ is finite. Let $%
\varphi :[0,\infty )\rightarrow \lbrack 0,\infty )$ be continuous strictly
increasing function. If $\ $ 
\begin{equation}
\varphi \left( a\right) \otimes c\geq \varphi \left( a\otimes c\right) ,
\label{poiu}
\end{equation}%
then the inequality 
\begin{equation*}
\mathbf{I}_{\otimes }\left( m,\varphi \left( f\right) \right) \geq \varphi
\left( \mathbf{I}_{\otimes }\left( m,f\right) \right)
\end{equation*}%
holds.
\end{corl}

Again, if we take $\varphi _{1}\left( x\right) =x$ in Theorem \ref{th3-3},
then we have the reverse Jensen inequality for universal integral.

\begin{corl}
Let $f\in $ $\mathcal{F}^{(X,\mathcal{A})}$ be a measurable function and $%
\otimes \colon \lbrack 0,\infty ]\sp{2}\rightarrow \lbrack 0,\infty ]$ be
the pseudo-multiplication with neutral element $e\in (0,\infty ]$ and $m$ $%
\in \mathcal{M}^{(X,\mathcal{A})}$ be a monotone measure such that $\mathbf{I%
}_{\otimes }\left( m,\varphi \left( f\right) \right) $ is finite. Let $%
\varphi :[0,\infty )\rightarrow \lbrack 0,\infty )$ be continuous strictly
increasing function. If $\ $ 
\begin{equation}
\varphi \left( a\otimes c\right) \geq \left( \varphi \left( a\right) \otimes
c\right) ,  \label{loil}
\end{equation}%
then the inequality 
\begin{equation*}
\varphi \left( \mathbf{I}_{\otimes }\left( m,f\right) \right) \geq \mathbf{I}%
_{\otimes }\left( m,\varphi \left( f\right) \right)
\end{equation*}%
holds.
\end{corl}

\begin{remk}
\label{remark3.6} If $\varphi :[0,\infty )\rightarrow \lbrack 0,\infty )$ is
continuous strictly increasing function such that $\varphi \left( x\right)
\leq x$ for all $x\in \lbrack 0,\infty )$ and $\varphi \;$is subdistributive
over $\ \otimes $, then (\ref{poiu}) holds readily. Indeed, 
\begin{equation*}
\varphi \left( a\otimes c\right) \leq \varphi \left( a\right) \otimes
\varphi \left( c\right) \leq \varphi \left( a\right) \otimes c.
\end{equation*}%
Also, if $\varphi \left( x\right) \geq x$ for all $x\in \lbrack 0,\infty )$
and $\varphi \;$is superdistributive over $\ \otimes $, then (\ref{loil})
holds similarly, i.e.,%
\begin{equation*}
\varphi \left( a\otimes c\right) \geq \varphi \left( a\right) \otimes
\varphi \left( c\right) \geq \varphi \left( a\right) \otimes c.
\end{equation*}
\end{remk}

\begin{corl}
Let $f\in $ $\mathcal{F}^{(X,\mathcal{A})}$ be a measurable function and $%
\otimes \colon \lbrack 0,\infty ]\sp{2}\rightarrow \lbrack 0,\infty ]$ be
the pseudo-multiplication with neutral element $e\in (0,\infty ]$ and $m$ $%
\in \mathcal{M}^{(X,\mathcal{A})}$ be a monotone measure such that $\mathbf{I%
}_{\otimes }\left( m,f\right) $ is finite. Let $\varphi :[0,\infty
)\rightarrow \lbrack 0,\infty )$ be continuous strictly increasing function
such that $\varphi \left( x\right) \leq x$ for all $x\in \lbrack 0,\infty )$%
. Then the inequality 
\begin{equation*}
\mathbf{I}_{\otimes }\left( m,\varphi \left( f\right) \right) \geq \varphi
\left( \mathbf{I}_{\otimes }\left( m,f\right) \right)
\end{equation*}%
holds, where $\varphi \;$is subdistributive over $\ \otimes $.
\end{corl}

\begin{corl}
Let $f\in $ $\mathcal{F}^{(X,\mathcal{A})}$ be a measurable function and $%
\otimes \colon \lbrack 0,\infty ]\sp{2}\rightarrow \lbrack 0,\infty ]$ be
the pseudo-multiplication with neutral element $e\in (0,\infty ]$ and $m$ $%
\in \mathcal{M}^{(X,\mathcal{A})}$ be a monotone measure such that $\mathbf{I%
}_{\otimes }\left( m,\varphi \left( f\right) \right) $ is finite. Let $%
\varphi :[0,\infty )\rightarrow \lbrack 0,\infty )$ be continuous strictly
increasing function such that $\varphi \left( x\right) \geq x$ for all $x\in
\lbrack 0,\infty )$. Then the inequality 
\begin{equation*}
\varphi \left( \mathbf{I}_{\otimes }\left( m,f\right) \right) \geq \mathbf{I}%
_{\otimes }\left( m,\varphi \left( f\right) \right)
\end{equation*}%
holds, where $\varphi \;$is superdistributive over $\ \otimes $.
\end{corl}

Notice that if the pseudo-multiplication $\otimes $ is minimum (i.e., for
Sugeno integral), then the following results hold (see \cite{RomFloCha07b}
for asimilar result).

\begin{corl}
\cite{RomFloCha07b} Let $f\in $ $\mathcal{F}^{(X,\mathcal{A})}$ be a
measurable function and $m$ $\in \mathcal{M}^{(X,\mathcal{A})}$ be a
monotone measure such that $\mathbf{Su}\left( m,f\right) $ is finite. Let $%
\varphi :[0,\infty )\rightarrow \lbrack 0,\infty )$ be continuous strictly
increasing function such that $\varphi \left( x\right) \leq x$ for all $x\in
\lbrack 0,\infty )$. Then the inequality 
\begin{equation*}
\mathbf{Su}\left( m,\varphi \left( f\right) \right) \geq \varphi \left( 
\mathbf{Su}\left( m,f\right) \right)
\end{equation*}%
holds.
\end{corl}

\begin{corl}
\cite{RomFloCha07b} Let $f\in $ $\mathcal{F}^{(X,\mathcal{A})}$ be a
measurable function and $m$ $\in \mathcal{M}^{(X,\mathcal{A})}$ be a
monotone measure such that $\mathbf{Su}\left( m,\varphi \left( f\right)
\right) $ is finite. Let $\varphi :[0,\infty )\rightarrow \lbrack 0,\infty )$
be continuous strictly increasing function such that $\varphi \left(
x\right) \geq x$ for all $x\in \lbrack 0,\infty )$. Then the inequality 
\begin{equation*}
\varphi \left( \mathbf{Su}\left( m,f\right) \right) \geq \mathbf{Su}\left(
m,\varphi \left( f\right) \right)
\end{equation*}%
holds.
\end{corl}

\bigskip When $\varphi _{1}\left( x\right) =x^{s}$ and $\varphi _{2}\left(
x\right) =x^{r}$ for all $r,s\in \left( 0,\infty \right) $ in Theorem \ref%
{th3-3}, then we have the following Lyapunov inequality for universal
integral.

\begin{corl}
Let $f\in $ $\mathcal{F}^{(X,\mathcal{A})}$ be a measurable function and $%
\otimes \colon \lbrack 0,\infty ]\sp{2}\rightarrow \lbrack 0,\infty ]$ be
the pseudo-multiplication with neutral element $e\in (0,\infty ]$ and $m$ $%
\in \mathcal{M}^{(X,\mathcal{A})}$ be a monotone measure such that $\mathbf{I%
}_{\otimes }\left( m,f^{r}\right) $ is finite. If $\ $ 
\begin{equation*}
\left( a^{s}\otimes c\right) ^{\frac{1}{s}}\geq \left( a^{r}\otimes c\right)
^{\frac{1}{r}},
\end{equation*}%
then the inequality 
\begin{equation*}
\left( \mathbf{I}_{\otimes }\left( m,f^{s}\right) \right) ^{\frac{1}{s}}\geq
\left( \mathbf{I}_{\otimes }\left( m,f^{r}\right) \right) ^{\frac{1}{r}}
\end{equation*}%
holds for all $r,s\in \left( 0,\infty \right) $.
\end{corl}

Notice that when working on $[0,1]$ in Theorem \ref{th3-3}, we mostly deal
with $e=1$, then $\otimes =\circledast $ is semicopula ($t$-seminorm) and
the following results hold.

\begin{corl}
Let $f\in $ $\mathcal{F}_{\left[ 0,1\right] }^{(X,\mathcal{A})}$ be a
measurable function and $m$ $\in \mathcal{M}_{1}^{(X,\mathcal{A})}$ be a
monotone measure. Let $\varphi _{i}:[0,\infty )\rightarrow \lbrack 0,\infty
),i=1,2$ be continuous strictly increasing functions. If semicopula $%
\circledast $ satisfies $\ $ 
\begin{equation*}
\varphi _{1}^{-1}\left( \varphi _{1}\left( a\right) \circledast c\right)
\geq \varphi _{2}^{-1}\left( \varphi _{2}\left( a\right) \circledast
c\right) ,
\end{equation*}%
then the inequality 
\begin{equation*}
\varphi _{1}^{-1}\left( \mathbf{I}_{\circledast }\left( m,\varphi _{1}\left(
f\right) \right) \right) \geq \varphi _{2}^{-1}\left( \mathbf{I}%
_{\circledast }\left( m,\varphi _{2}\left( f\right) \right) \right)
\end{equation*}%
holds.
\end{corl}

\begin{corl}
\label{cor3.16} Let $f\in $ $\mathcal{F}_{\left[ 0,1\right] }^{(X,\mathcal{A}%
)}$ be a measurable function and $m$ $\in \mathcal{M}_{1}^{(X,\mathcal{A})}$
be a monotone measure. Let $\varphi :[0,1]\rightarrow \lbrack 0,1]$ be
continuous strictly increasing function such that $\varphi \left( x\right)
\leq x$ for all $x\in \lbrack 0,1]$. Then the inequality 
\begin{equation*}
\mathbf{I}_{\circledast }\left( m,\varphi \left( f\right) \right) \geq
\varphi \left( \mathbf{I}_{\circledast }\left( m,f\right) \right)
\end{equation*}%
holds, where $\varphi \;$is subdistributive over semicopula $\circledast $.
\end{corl}

\begin{corl}
\label{cor3.15} Let $f\in $ $\mathcal{F}_{\left[ 0,1\right] }^{(X,\mathcal{A}%
)}$ be a measurable function and $m$ $\in \mathcal{M}_{1}^{(X,\mathcal{A})}$
be a monotone measure. Let $\varphi :[0,1]\rightarrow \lbrack 0,1]$ be
continuous strictly increasing function such that $\varphi \left( x\right)
\geq x$ for all $x\in \lbrack 0,1]$. Then the inequality 
\begin{equation*}
\varphi \left( \mathbf{I}_{\circledast }\left( m,f\right) \right) \geq 
\mathbf{I}_{\circledast }\left( m,\varphi \left( f\right) \right)
\end{equation*}%
holds, where $\varphi \;$is superdistributive over semicopula $\circledast $.
\end{corl}

\begin{corl}
Let $f\in $ $\mathcal{F}_{\left[ 0,1\right] }^{(X,\mathcal{A})}$ be a
measurable function and $m$ $\in \mathcal{M}_{1}^{(X,\mathcal{A})}$ be a
monotone measure. If semicopula $\circledast $ satisfies $\ $ 
\begin{equation}
\left( a^{s}\circledast c\right) ^{\frac{1}{s}}\geq \left( a^{r}\circledast
c\right) ^{\frac{1}{r}},
\end{equation}%
then the inequality 
\begin{equation*}
\left( \mathbf{I}_{\circledast }\left( m,f^{s}\right) \right) ^{\frac{1}{s}%
}\geq \left( \mathbf{I}_{\circledast }\left( m,f^{r}\right) \right) ^{\frac{1%
}{r}}
\end{equation*}%
holds for all $r,s\in \left( 0,\infty \right) $.
\end{corl}

\begin{corl}
Let $f\in $ $\mathcal{F}_{\left[ 0,1\right] }^{(X,\mathcal{A})}$ be a
measurable function and $m$ $\in \mathcal{M}_{1}^{(X,\mathcal{A})}$ be a
monotone measure. then the inequality 
\begin{equation*}
\left( \mathbf{Su}\left( m,f^{s}\right) \right) ^{\frac{1}{s}}\geq \left( 
\mathbf{Su}\left( m,f^{r}\right) \right) ^{\frac{1}{r}}
\end{equation*}%
holds for all $0<r\leq s<\infty $.
\end{corl}

\section{On reverse inequalities}

\hspace{0.5cm}By using the concepts of t-seminorm and t-semiconorm, Su\'{a}%
rez and Gil proposed the a family of semiconormed integrals \cite{SurGil}.
Define 
\begin{equation*}
\mathbf{I}_{\oplus }\left( m,f\right) =\inf \left\{ t\oplus m\left( \left\{
f>t\right\} \right) \;|\;t\in (0,\infty ])\right\} .
\end{equation*}%
Hence, we get the following theorems.

\begin{thrm}
\label{thph}Let a non-decreasing $n$-place function $H:[0,\infty
)^{n}\rightarrow \lbrack 0,\infty )$ such that $H$ be continuous. If $\
\oplus \colon \lbrack 0,\infty ]\sp{n}\rightarrow \lbrack 0,\infty ]$ is the
pseudo-addition with neutral element $0$, satisfies 
\begin{eqnarray*}
&&U_{0}^{-1}\left[ U_{0}\left( H\left( \psi _{1}\left( p_{1}\right) ,\psi
_{2}\left( p_{2}\right) ,...,\psi _{n}\left( p_{n}\right) \right) \right)
\oplus c\right] \overset{}{\leq }H\left( \psi _{1}\left( U_{1}^{-1}\left[
\left( U_{1}\left( p_{1}\right) \right) \oplus c\right] \right) ,\psi
_{2}\left( p_{2}\right) ,...,\psi _{n}\left( p_{n}\right) \right) \\
&&\wedge H\left( \psi _{1}\left( p_{1}\right) ,\psi _{2}\left( U_{2}^{-1} 
\left[ \left( U_{2}\left( p_{2}\right) \right) \oplus c\right] \right) ,\psi
_{3}\left( p_{3}\right) ,...,\psi _{n}\left( p_{n}\right) \right) \\
&&\wedge ...\wedge H\left( \psi _{1}\left( p_{1}\right) ,\psi _{2}\left(
p_{2}\right) ,...,\psi _{n-1}\left( p_{n-1}\right) ,\psi \left( U_{n}^{-1} 
\left[ \left( U_{n}\left( p_{n}\right) \right) \oplus c\right] \right)
\right) ,
\end{eqnarray*}%
then for any system $U_{0},U_{1},...,U_{n}:[0,\infty )\rightarrow \lbrack
0,\infty )$ of continuous strictly increasing functions, and any system $%
\psi _{1},\psi _{2},...,\psi _{n}:[0,\infty )\rightarrow \lbrack 0,\infty )$
of continuous increasing functions and any comontone system $%
f_{1},f_{2},...,f_{n}\in $ $\mathcal{F}^{(X,\mathcal{A})}$ and a monotone
measure $m$ $\in \mathcal{M}^{(X,\mathcal{A})}$, $\mathbf{I}_{\oplus }\left(
m,U_{i}\left( f_{i}\right) \right) <\infty $ for all $i=1,2,...n$, it holds%
\begin{equation*}
U_{0}^{-1}[\mathbf{I}_{\oplus }\left( m,U_{0}[H\left( \psi _{1}\left(
f_{1}\right) ,...,\psi _{n}\left( f_{n}\right) \right) ]\right) ]\overset{}{%
\leq }H\left[ \psi _{1}\left( U_{1}^{-1}\left( \mathbf{I}_{\oplus }\left(
m,U_{1}\left( f_{1}\right) \right) \right) \right) ,...,\psi _{n}\left(
U_{n}^{-1}\left( \mathbf{I}_{\oplus }\left( m,U_{n}\left( f_{n}\right)
\right) \right) \right) \right] .
\end{equation*}
\end{thrm}

\noindent \textbf{Proof.}\ Let $\mathbf{I}_{\oplus }\left( m,U_{i}\left(
f_{i}\right) \right) =p_{i}<\infty $ for all $i=1,2,...,n$. So, for any $%
\varepsilon >0$, there exist $p_{i(\varepsilon )}$ such that 
\begin{equation*}
m(\left\{ U_{i}\left( f_{i}\right) >p_{i(\varepsilon )}\right\} )=M\sb{i},
\end{equation*}%
where $p_{i(\varepsilon )}\oplus M\sb{i}\leq p_{i}+\varepsilon $ for all $%
i=1,2,...,n.$ Then,%
\begin{equation*}
\psi _{i}\left( U_{i}^{-1}\left[ p_{i\left( \varepsilon \right) }\oplus M%
\sb{i}\right] \right) \leq \psi _{i}\left( U_{i}^{-1}\left[
p_{i}+\varepsilon \right] \right) ,\text{ for all }i=1,2,...,n.
\end{equation*}%
Then,%
\begin{equation*}
\psi _{i}\left( U_{i}^{-1}\left[ p_{i\left( \varepsilon \right) }\right]
\right) =\psi _{i}\left( U_{i}^{-1}\left[ p_{i\left( \varepsilon \right)
}\oplus 0\right] \right) \leq \psi _{i}\left( U_{i}^{-1}\left[
p_{i}+\varepsilon \right] \right) ,\text{ for all }i=1,2,...,n.
\end{equation*}%
The comonotonicity of $f_{1},f_{2},...,f_{n}$ and the monotonicity of $H$
imply that 
\begin{eqnarray*}
&&m\left( \{U_{0}\left( H\left( \psi _{1}\left( f_{1}\right) ,...,\psi
_{n}\left( f_{n}\right) \right) \right) >U_{0}\left( H\left( \psi _{1}\left(
U_{1}^{-1}\left( p_{1(\varepsilon )}\right) \right) ,...,\psi _{n}\left(
U_{n}^{-1}\left( p_{n(\varepsilon )}\right) \right) \right) \right) \}\right)
\\
&=&m(\{H\left( \psi _{1}\left( f_{1}\right) ,...,\psi _{n}\left(
f_{n}\right) \right) >H\left( \psi _{1}\left( U_{1}^{-1}\left(
p_{1(\varepsilon )}\right) \right) ,...,\psi _{n}\left( U_{n}^{-1}\left(
p_{n(\varepsilon )}\right) \right) \right) \}) \\
&\leq &m(\left\{ U_{1}\left( f_{1}\right) >p_{1(\varepsilon )}\right\} )\vee
m(\left\{ U_{2}\left( f_{2}\right) >p_{2(\varepsilon )}\right\} )\vee
....\vee m(\left\{ U_{n}\left( f_{n}\right) >p_{n(\varepsilon )}\right\} ) \\
&=&M_{1}\vee M_{2}\vee ...\vee M_{n}.
\end{eqnarray*}%
Hence 
\begin{eqnarray*}
&&U_{0}^{-1}\left[ \inf \left( t\oplus m(\{U_{0}\left( H\left( \psi
_{1}\left( f_{1}\right) ,...,\psi _{n}\left( f_{n}\right) \right) \right)
>t\})\;|\;t\in (0,\infty ])\right) \right] \\
&&\overset{}{\leq }U_{0}^{-1}\left( \left[ 
\begin{array}{c}
U_{0}\left( H\left( \psi _{1}\left( U_{1}^{-1}\left( p_{1(\varepsilon
)}\right) \right) ,...,\psi _{n}\left( U_{n}^{-1}\left( p_{n(\varepsilon
)}\right) \right) \right) \right) \oplus \\ 
m(\{U_{0}\left( H\left( \psi _{1}\left( f_{1}\right) ,...,\psi _{n}\left(
f_{n}\right) \right) \right) >U_{0}\left( H\left( \psi _{1}\left(
U_{1}^{-1}\left( p_{1(\varepsilon )}\right) \right) ,...,\psi _{n}\left(
U_{n}^{-1}\left( p_{n(\varepsilon )}\right) \right) \right) \right) \})%
\end{array}%
\right] \right) \\
&&\overset{}{\leq }U_{0}^{-1}\left( \left[ U_{0}\left( H\left( \psi
_{1}\left( U_{1}^{-1}\left( p_{1(\varepsilon )}\right) \right) ,...,\psi
_{n}\left( U_{n}^{-1}\left( p_{n(\varepsilon )}\right) \right) \right)
\right) \oplus \left( M_{1}\vee M_{2}\vee ...\vee M_{n}\right) \right]
\right) \\
&&\overset{}{=}\left( 
\begin{array}{c}
U_{0}^{-1}\left[ U_{0}\left( H\left( \psi _{1}\left( U_{1}^{-1}\left(
p_{1(\varepsilon )}\right) \right) ,...,\psi _{n}\left( U_{n}^{-1}\left(
p_{n(\varepsilon )}\right) \right) \right) \right) \oplus M_{1}\right] \\ 
\vee U_{0}^{-1}\left[ U_{0}\left( H\left( \psi _{1}\left( U_{1}^{-1}\left(
p_{1(\varepsilon )}\right) \right) ,...,\psi _{n}\left( U_{n}^{-1}\left(
p_{n(\varepsilon )}\right) \right) \right) \right) \oplus M_{2}\right] \\ 
\vee ...\vee U_{0}^{-1}\left[ U_{0}\left( H\left( \psi _{1}\left(
U_{1}^{-1}\left( p_{1(\varepsilon )}\right) \right) ,...,\psi _{n}\left(
U_{n}^{-1}\left( p_{n(\varepsilon )}\right) \right) \right) \right) \oplus
M_{n}\right]%
\end{array}%
\right) \\
&&\overset{}{\leq }\left( 
\begin{array}{c}
H\left( \psi _{1}\left( U_{1}^{-1}\left[ p_{1(\varepsilon )}\oplus M_{1}%
\right] \right) ,\psi _{2}\left( U_{2}^{-1}\left( p_{2(\varepsilon )}\right)
\right) ,...,\psi _{n}\left( U_{n}^{-1}\left( p_{n(\varepsilon )}\right)
\right) \right) \\ 
\vee H\left( \psi _{1}\left( U_{1}^{-1}\left( p_{1(\varepsilon )}\right)
\right) ,\psi _{2}\left( U_{2}^{-1}\left[ p_{2(\varepsilon )}\oplus M_{2}%
\right] \right) ,\psi _{3}\left( U_{3}^{-1}\left( p_{3(\varepsilon )}\right)
\right) ,...,\psi _{n}\left( U_{n}^{-1}\left( p_{n(\varepsilon )}\right)
\right) \right) \\ 
\vee ...\vee H\left( \psi _{1}\left( U_{1}^{-1}\left( p_{1(\varepsilon
)}\right) \right) ,...,\psi _{n-1}\left( U_{n-1}^{-1}\left(
p_{(n-1)(\varepsilon )}\right) \right) ,\psi _{n}\left( U_{n}^{-1}\left[
p_{n(\varepsilon )}\oplus M_{n}\right] \right) \right)%
\end{array}%
\right) \\
&&\overset{}{\leq }\left( 
\begin{array}{c}
H\left( \psi _{1}\left( U_{1}^{-1}\left[ p_{1}+\varepsilon \right] \right)
,\psi _{2}\left( U_{2}^{-1}\left( p_{2(\varepsilon )}\right) \right)
,...,\psi _{n}\left( U_{n}^{-1}\left( p_{n(\varepsilon )}\right) \right)
\right) \\ 
\vee H\left( \psi _{1}\left( U_{1}^{-1}\left( p_{1(\varepsilon )}\right)
\right) ,\psi _{2}\left( U_{2}^{-1}\left[ p_{2}+\varepsilon \right] \right)
,\psi _{3}\left( U_{3}^{-1}\left( p_{3(\varepsilon )}\right) \right)
,...,\psi _{n}\left( U_{n}^{-1}\left( p_{n(\varepsilon )}\right) \right)
\right) \\ 
\vee ...\vee H\left( \psi _{1}\left( U_{1}^{-1}\left( p_{1(\varepsilon
)}\right) \right) ,...,\psi _{n-1}\left( U_{n-1}^{-1}\left(
p_{(n-1)(\varepsilon )}\right) \right) ,\psi _{n}\left( U_{n}^{-1}\left[
p_{n}+\varepsilon \right] \right) \right)%
\end{array}%
\right) \\
&&\overset{}{\leq }H\left( \psi _{1}\left( U_{1}^{-1}\left[
p_{1}+\varepsilon \right] \right) ,\psi _{2}\left( U_{2}^{-1}\left[
p_{2}+\varepsilon \right] \right) ,...,\psi _{n}\left( U_{n}^{-1}\left[
p_{n}+\varepsilon \right] \right) \right) ,
\end{eqnarray*}%
whence $U_{0}^{-1}[\mathbf{I}_{\oplus }\left( m,U_{0}[H\left( \psi
_{1}\left( f_{1}\right) ,...,\psi _{n}\left( f_{n}\right) \right) ]\right)
]\leq H\left( \psi _{1}\left( U_{1}^{-1}\left[ p_{1}\right] \right) ,\psi
_{2}\left( U_{2}^{-1}\left[ p_{2}\right] \right) ,...,\psi _{n}\left(
U_{n}^{-1}\left[ p_{n}\right] \right) \right) $ follows from the continuity
of $H,\psi _{i},U_{i}$ for all $i,$ and the arbitrariness of $\varepsilon $.
And the theorem is proved. $\Box $

\begin{corl}
Let a non-decreasing $n$-place function $H:[0,1]^{n}\rightarrow \lbrack 0,1]$
such that $H$ be continuous and a continuous non-decreasing $\psi
:[0,1]\rightarrow \lbrack 0,1]$ be given. If the $t$-semiconorm $S$
satisfies 
\begin{eqnarray*}
&&U_{0}^{-1}\left[ S\left( U_{0}\left( H\left( \psi \left( p_{1}\right)
,\psi \left( p_{2}\right) ,...,\psi \left( p_{n}\right) \right) \right)
,c\right) \right] \overset{}{\leq }H\left( \psi \left( U_{1}^{-1}\left[
S\left( U_{1}\left( p_{1}\right) ,c\right) \right] \right) ,\psi \left(
p_{2}\right) ,...,\psi \left( p_{n}\right) \right) \\
&&\wedge H\left( \psi \left( p_{1}\right) ,\psi \left( U_{2}^{-1}\left[
S\left( U_{2}\left( p_{2}\right) ,c\right) \right] \right) ,\psi \left(
p_{3}\right) ,...,\psi \left( p_{n}\right) \right) \\
&&\wedge ...\wedge H\left( \psi \left( p_{1}\right) ,\psi \left(
p_{2}\right) ,...,\psi \left( p_{\left( n-1\right) }\right) ,\psi \left(
U_{n}^{-1}\left[ S\left( U_{n}\left( p_{n}\right) ,c\right) \right] \right)
\right) ,
\end{eqnarray*}%
then for any system $U_{0},U_{1},...,U_{n}:[0,1]\rightarrow \lbrack 0,1]$ of
continuous strictly increasing functions and any comontone system $%
f_{1},f_{2},...,f_{n}\in $ $\mathcal{F}_{\left[ 0,1\right] }^{(X,\mathcal{A}%
)}$ and a monotone measure $m\in \mathcal{M}_{1}^{(X,\mathcal{A})}$, it holds%
\begin{equation*}
U_{0}^{-1}[\mathbf{I}_{S}\left( m,U_{0}[H\left( \psi \left( f_{1}\right)
,...,\psi \left( f_{n}\right) \right) ]\right) ]\leq H\left[ \psi \left(
U_{1}^{-1}\left( \mathbf{I}_{S}\left( m,U_{1}\left( f_{1}\right) \right)
\right) \right) ,...,\psi \left( U_{n}^{-1}\left( \mathbf{I}_{S}\left(
m,U_{n}\left( f_{n}\right) \right) \right) \right) \right] .
\end{equation*}
\end{corl}

In an analogous way as in the proof of Theorem \ref{thph} we have the
following results.

\begin{thrm}
Let a non-decreasing $n$-place function $H:[0,\infty )^{n}\rightarrow
\lbrack 0,\infty )$ such that $H$ be continuous. If $\ \oplus \colon \lbrack
0,\infty ]\sp{n}\rightarrow \lbrack 0,\infty ]$ is the pseudo-addition with
neutral element $0$, satisfies 
\begin{eqnarray*}
&&\left( \left( H\left( p_{1},p_{2},...,p_{n}\right) \right) ^{\xi
_{0}}\oplus c\right) ^{\omega _{0}}\overset{}{\leq }H\left( \left(
p_{1}^{\xi _{1}}\oplus c\right) ^{\omega _{1}},p_{2},...,p_{n}\right) \wedge
\\
&&H\left( p_{1},\left( p_{2}^{\xi _{2}}\oplus c\right) ^{\omega
_{2}},p_{3},...,p_{n}\right) \wedge ...\wedge H\left(
p_{1},p_{2},...,p_{n-1},\left( p_{n}^{\xi _{n}}\oplus c\right) ^{\omega
_{n}}\right) ,
\end{eqnarray*}%
then for any comontone system $f_{1},f_{2},...,f_{n}\in $ $\mathcal{F}^{(X,%
\mathcal{A})}$ and a monotone measure $m$ $\in \mathcal{M}^{(X,\mathcal{A})}$%
, it holds%
\begin{equation*}
\left[ \mathbf{I}_{\oplus }\left( m,\left( H\left( f_{1},...,f_{n}\right)
\right) ^{\xi _{0}}\right) \right] ^{^{\omega _{0}}}\overset{}{\leq }H\left[
\left( \mathbf{I}_{\oplus }\left( m,f_{1}^{\xi _{1}}\right) \right) ^{\omega
_{1}},\left( \mathbf{I}_{\oplus }\left( m,f_{2}^{\xi _{2}}\right) \right)
^{\omega _{2}},...,\left( \mathbf{I}_{\oplus }\left( m,f_{n}^{\xi
_{n}}\right) \right) ^{\omega _{n}}\right]
\end{equation*}%
for all $\omega _{j},\xi _{j}\in \left( 0,\infty \right) ,$ $\omega _{i}\xi
_{i}\leq 1,$ where $i=1,2,...n$ and $j=0,1,2,...n$.
\end{thrm}

\begin{corl}
\label{th4-1} Let a non-decreasing $n$-place function $H:[0,\infty
)^{n}\rightarrow \lbrack 0,\infty )$ such that $H$ be continuous. If the $t$%
-semiconorm $S$ satisfies 
\begin{eqnarray*}
&&S^{\omega _{0}}\left( \left( H\left( p_{1},p_{2},...,p_{n}\right) \right)
^{\xi _{0}},c\right) \overset{}{\leq }H\left( S^{\omega _{1}}\left(
p_{1}^{\xi _{1}},c\right) ,p_{2},...,p_{n}\right) \wedge \\
&&H\left( p_{1},S^{\omega _{2}}\left( p_{2}^{\xi _{2}},c\right)
,p_{3},...,p_{n}\right) \wedge ...\wedge H\left(
p_{1},p_{2},...,p_{n-1},S^{\omega _{n}}\left( p_{n}^{\xi _{n}},c\right)
\right) ,
\end{eqnarray*}%
then for any comontone system $f_{1},f_{2},...,f_{n}\in $ $\mathcal{F}_{%
\left[ 0,1\right] }^{(X,\mathcal{A})}$ and a monotone measure $m$ $\in 
\mathcal{M}_{1}^{(X,\mathcal{A})}$, it holds%
\begin{equation*}
\left[ \mathbf{I}_{S}\left( m,\left( H\left( f_{1},...,f_{n}\right) \right)
^{\xi _{0}}\right) \right] ^{^{\omega _{0}}}\overset{}{\leq }H\left[ \left( 
\mathbf{I}_{S}\left( m,f_{1}^{\xi _{1}}\right) \right) ^{\omega _{1}},\left( 
\mathbf{I}_{S}\left( m,f_{2}^{\xi _{2}}\right) \right) ^{\omega
_{2}},...,\left( \mathbf{I}_{S}\left( m,f_{n}^{\xi _{n}}\right) \right)
^{\omega _{n}}\right]
\end{equation*}%
for all $\omega _{j},\xi _{j}\in \left( 0,\infty \right) ,$ $\omega _{i}\xi
_{i}\leq 1,$ where $i=1,2,...n$ and $j=0,1,2,...n$.
\end{corl}

\begin{corl}
Let $f,g\in $ $\mathcal{F}_{\left[ 0,1\right] }^{(X,\mathcal{A})}$ be two
comonotone measurable functions. Let $\star \colon \lbrack 0,1]\sp{2}%
\rightarrow \lbrack 0,1]$ be continuous and nondecreasing in both arguments.
If the semiconorm $S$ satisfies 
\begin{equation}
S^{\lambda }(\left( a\star b\right) ^{\alpha },c)\leq \left[ S^{\upsilon
}\left( a^{\beta },c\right) \star b\right] \wedge \left[ a\star S^{\tau
}(b^{\gamma },c)\right] ,  \label{0909}
\end{equation}%
then the inequality 
\begin{equation*}
\left[ \mathbf{I}_{S}\left( m,(f\star g)^{\alpha }\right) \right] ^{\lambda
}\leq \left[ \mathbf{I}_{S}\left( m,f^{\beta }\right) \right] ^{\upsilon
}\star \left[ \mathbf{I}_{S}\left( m,g^{\gamma }\right) \right] ^{\tau }
\end{equation*}%
holds for all $\alpha ,\beta ,\gamma ,\lambda ,\upsilon ,\tau \in \left(
0,\infty \right) ,\gamma \tau \leq 1,\beta \upsilon \leq 1$ and for any $%
m\in \mathcal{M}_{1}^{(X,\mathcal{A})}$.
\end{corl}

Let $\alpha =\beta =\gamma =k$ and $\lambda =\upsilon =\tau =\frac{1}{k}$
for all $k\in \left( 0,\infty \right) $, then we get the Minkowski
inequality for semiconormed fuzzy integrals (if $k=1,$ then we have the
reverse Chebyshev inequality for semiconormed fuzzy integrals \cite{OuyMes09}%
).

\begin{corl}
Let $f,g\in $ $\mathcal{F}_{\left[ 0,1\right] }^{(X,\mathcal{A})}$ be two
comonotone measurable functions. Let $\star \colon \lbrack 0,1]\sp{2}%
\rightarrow \lbrack 0,1]$ be continuous and nondecreasing in both arguments.
If the semiconorm $S$ satisfies 
\begin{equation*}
\left[ S(\left( a\star b\right) ^{k},c)\right] ^{\frac{1}{k}}\leq \left[
\left( S\left( a^{k},c\right) \right) ^{\frac{1}{k}}\star b\right] \wedge %
\left[ a\star \left( S(b^{k},c)\right) ^{\frac{1}{k}}\right] ,
\end{equation*}%
then the inequality 
\begin{equation*}
\left( \mathbf{I}_{S}\left( m,\left( f\star g\right) ^{k}\right) \right) ^{%
\frac{1}{k}}\leq \left( \mathbf{I}_{S}\left( m,f^{k}\right) \right) ^{\frac{1%
}{k}}\star \left( \mathbf{I}_{S}\left( m,g^{k}\right) \right) ^{\frac{1}{k}}
\end{equation*}%
holds for any $m\in \mathcal{M}_{1}^{(X,\mathcal{A})}$ and for all $%
0<k<\infty $.
\end{corl}

Notice that if the semiconorm $S$ is maximum (i.e., for Sugeno integral) and 
$\star $ is bounded from below by maximum, then $S$ is dominated by $\star $%
. Thus the following results hold.

\begin{corl}
Let $f,g\in $ $\mathcal{F}_{\left[ 0,1\right] }^{(X,\mathcal{A})}$ be two
comonotone measurable functions. Let $\star \colon \lbrack 0,1]\sp{2}%
\rightarrow \lbrack 0,1]$ be continuous and nondecreasing in both arguments
and bounded from below by maximum. Then the inequality 
\begin{equation*}
\left[ \mathbf{Su}\left( m,\left( f\star g\right) ^{\alpha }\right) \right]
^{\lambda }\leq \left[ \mathbf{Su}\left( m,f^{\beta }\right) \right]
^{\upsilon }\star \left[ \mathbf{Su}\left( m,g^{\gamma }\right) \right]
^{\tau }
\end{equation*}%
holds for any $m\in \mathcal{M}_{1}^{(X,\mathcal{A})}$ and for all $\alpha
,\beta ,\gamma ,\lambda ,\upsilon ,\tau \in \left( 0,\infty \right) ,1\leq
\alpha \lambda <\infty ,0<\beta \upsilon \leq 1,0<\gamma \tau \leq 1,\lambda
\geq \tau ,\upsilon .$
\end{corl}

\begin{corl}
(\textrm{\cite{Aga09}) }Let $f,g\in $ $\mathcal{F}_{\left[ 0,1\right] }^{(X,%
\mathcal{A})}$ be two comonotone measurable functions. Let $\star \colon
\lbrack 0,1]\sp{2}\rightarrow \lbrack 0,1]$ be continuous and nondecreasing
in both arguments and bounded from below by maximum. Then the inequality 
\begin{equation*}
\left( \mathbf{Su}\left( m,\left( f\star g\right) ^{k}\right) \right) ^{%
\frac{1}{k}}\leq \left( \mathbf{Su}\left( m,f^{k}\right) \right) ^{\frac{1}{k%
}}\star \left( \mathbf{Su}\left( m,g^{k}\right) \right) ^{\frac{1}{k}}
\end{equation*}%
holds for any $m\in \mathcal{M}_{1}^{(X,\mathcal{A})}$ and for all $%
0<k<\infty $.
\end{corl}

\begin{corl}
Let $f,g\in $ $\mathcal{F}_{\left[ 0,1\right] }^{(X,\mathcal{A})}$ be two
comonotone measurable functions. Let $\star \colon \lbrack 0,1]\sp{2}%
\rightarrow \lbrack 0,1]$ be continuous and nondecreasing in both arguments
and bounded from below by maximum. Then the inequality 
\begin{equation*}
\mathbf{Su}\left( m,\left( f\star g\right) \right) \leq \left( \mathbf{Su}%
\left( m,f^{p}\right) \right) ^{\frac{1}{p}}\star \left( \mathbf{Su}\left(
m,g^{q}\right) \right) ^{\frac{1}{q}}
\end{equation*}%
holds for any $m\in \mathcal{M}_{1}^{(X,\mathcal{A})}$ and for all $p,q\in
\lbrack 1,\infty )$.
\end{corl}

\begin{corl}
Let $f,g\in $ $\mathcal{F}_{\left[ 0,1\right] }^{(X,\mathcal{A})}$ be two
comonotone measurable functions. Let $\star \colon \lbrack 0,1]\sp{2}%
\rightarrow \lbrack 0,1]$ be continuous and nondecreasing in both arguments
and bounded from below by maximum. Then the inequality 
\begin{equation*}
\mathbf{Su}\left( m,\left( f\star g\right) \right) \leq \mathbf{Su}\left(
m,f\right) \star \mathbf{Su}\left( m,g\right)
\end{equation*}%
holds for any $m\in \mathcal{M}_{1}^{(X,\mathcal{A})}.$
\end{corl}

\begin{remk}
If $(x\star 0)\vee (0\star x)\geq x$ for any $x\in \lbrack 0,1]\ $and if $%
\Phi \left( x\right) =\left( .\right) ^{\alpha }$ is subdistributive over $%
\star $ and $S^{\lambda }$ dominates $\star $, then (\ref{0909}) holds
readily for all $\alpha ,\beta ,\gamma ,\lambda ,\upsilon ,\tau \in \left(
0,\infty \right) ,1\leq \alpha \lambda <\infty ,0<\beta \upsilon \leq
1,0<\gamma \tau \leq 1,\alpha \geq \beta ,\gamma $ and $\lambda \geq \tau
,\upsilon .$
\end{remk}

Suppose the semiconorm $S$ further satisfies monotonicity and associativity
(i.e., it is a $t$-conorm). Then, we have the following result:

\begin{corl}
Let $(X,\mathcal{F},\mu )$ be a fuzzy measure space and $f,g\colon
X\rightarrow \lbrack 0,1]$ two comonotone measurable functions. If $S$ be a
continuous $t$-conorm, then 
\begin{equation*}
\left[ \mathbf{I}_{S}\left( m,S^{\alpha }\left( f,g\right) \right) \right]
^{\lambda }\leq S\left( \left[ \mathbf{I}_{S}\left( m,f^{\beta }\right) %
\right] ^{\upsilon },\left[ \mathbf{I}_{S}\left( m,g^{\gamma }\right) \right]
^{\tau }\right)
\end{equation*}%
holds for any $m\in \mathcal{M}_{1}^{(X,\mathcal{A})}$ and for all $\alpha
,\beta ,\gamma ,\lambda ,\upsilon ,\tau \in \left( 0,\infty \right) ,1\leq
\alpha \lambda <\infty ,0<\beta \upsilon \leq 1,0<\gamma \tau \leq 1,\alpha
\geq \beta ,\gamma $ and $\lambda \geq \tau ,\upsilon ,$ where $\left(
.\right) ^{\alpha }$ is subdistributive over $S$, $S^{\lambda }$ dominates $%
S $ and $S(f,g)(x)=S(f(x),g(x))$ for any $x\in X$.
\end{corl}

Let $\alpha =\beta =\gamma =\lambda =\upsilon =\tau =1$, then we have the
following result:

\begin{corl}
Let $(X,\mathcal{F},\mu )$ be a fuzzy measure space and $f,g\colon
X\rightarrow \lbrack 0,1]$ two comonotone measurable functions. If $S$ be a
continuous $t$-conorm, then 
\begin{equation*}
\mathbf{I}_{S}\left( m,S\left( f,g\right) \right) \leq S\left( \mathbf{I}%
_{S}\left( m,f\right) ,\mathbf{I}_{S}\left( m,g\right) \right)
\end{equation*}%
holds for any $m\in \mathcal{M}_{1}^{(X,\mathcal{A})},$ where $%
S(f,g)(x)=S(f(x),g(x))$ for any $x\in X$.
\end{corl}

\begin{thrm}
\label{th3-3}Let $f\in $ $\mathcal{F}^{(X,\mathcal{A})}$ be a measurable
function and $\oplus \colon \lbrack 0,\infty ]\sp{n}\rightarrow \lbrack
0,\infty ]$ be the pseudo-addition with neutral element $0$, satisfies and $%
m $ $\in \mathcal{M}^{(X,\mathcal{A})}$ be a monotone measure such that $%
\mathbf{I}_{\oplus }\left( m,\varphi _{1}\left( f\right) \right) $ is
finite. Let $\varphi _{i}:[0,\infty )\rightarrow \lbrack 0,\infty ),i=1,2$
be continuous strictly increasing functions. If $\ $ 
\begin{equation*}
\varphi _{1}^{-1}\left( \varphi _{1}\left( a\right) \oplus c\right) \leq
\varphi _{2}^{-1}\left( \varphi _{2}\left( a\right) \oplus c\right) ,
\end{equation*}%
then the inequality 
\begin{equation*}
\varphi _{1}^{-1}\left( \mathbf{I}_{\oplus }\left( m,\varphi _{1}\left(
f\right) \right) \right) \leq \varphi _{2}^{-1}\left( \mathbf{I}_{\oplus
}\left( m,\varphi _{2}\left( f\right) \right) \right)
\end{equation*}%
holds.
\end{thrm}

\noindent\textbf{Proof.}\ Let $\mathbf{I}_{\oplus }\left( m,\varphi _{2}\left(
f\right) \right) =p<\infty $. So, for any $\varepsilon >0$, there exists $%
p_{\varepsilon }$ such that $m(\left\{ \varphi _{2}\left( f\right) \geq
p_{\varepsilon }\right\} )=M,$where $p_{\varepsilon }\oplus M\leq
p+\varepsilon $. Hence, 
\begin{eqnarray*}
&&\varphi _{1}^{-1}\left( \mathbf{I}_{\otimes }\left( m,\varphi _{1}\left(
f\right) \right) \right) \overset{}{\leq }\varphi _{1}^{-1}\left( \left[
\varphi _{1}\left( \varphi _{2}^{-1}\left( p_{\varepsilon }\right) \right)
\oplus m(\{\varphi _{1}\left( f\right) \geq \varphi _{1}\left( \varphi
_{2}^{-1}\left( p_{\varepsilon }\right) \right) \})\right] \right) \\
&&\overset{}{=}\varphi _{1}^{-1}\left( \left[ \varphi _{1}\left( \varphi
_{2}^{-1}\left( p_{\varepsilon }\right) \right) \oplus m(\{\varphi
_{2}\left( f\right) \geq p_{\varepsilon }\})\right] \right) \\
&&\overset{}{\leq }\varphi _{2}^{-1}\left( \left[ \varphi _{2}\left( \varphi
_{2}^{-1}\left( p_{\varepsilon }\right) \right) \oplus m(\{\varphi
_{2}\left( f\right) \geq p_{\varepsilon }\})\right] \right) \\
&&\overset{}{=}\varphi _{2}^{-1}\left( \left[ p_{\varepsilon }\oplus M\right]
\right) \overset{}{\leq }\varphi _{2}^{-1}\left( p+\varepsilon \right)
\end{eqnarray*}%
whence $\varphi _{1}^{-1}\left( \mathbf{I}_{\otimes }\left( m,\varphi
_{1}\left( f\right) \right) \right) \geq \varphi _{2}^{-1}\left( p\right) $
follows from the continuity of $\varphi _{2}$ and the arbitrariness of $%
\varepsilon $. And the theorem is proved. $\Box $

\begin{corl}
Let $f\in $ $\mathcal{F}_{\left[ 0,1\right] }^{(X,\mathcal{A})}$ be a
measurable function and $m$ $\in \mathcal{M}_{1}^{(X,\mathcal{A})}$ be a
monotone measure. Let $\varphi _{i}:[0,1]\rightarrow \lbrack 0,1],i=1,2$ be
continuous strictly increasing functions. If the semiconorm $S$ satisfiesIf $%
\ $ 
\begin{equation*}
\varphi _{1}^{-1}\left( S\left( \varphi _{1}\left( a\right) ,c\right)
\right) \leq \varphi _{2}^{-1}\left( S\left( \varphi _{2}\left( a\right)
,c\right) \right) ,
\end{equation*}%
then the inequality 
\begin{equation*}
\varphi _{1}^{-1}\left( \mathbf{I}_{S}\left( m,\varphi _{1}\left( f\right)
\right) \right) \leq \varphi _{2}^{-1}\left( \mathbf{I}_{S}\left( m,\varphi
_{2}\left( f\right) \right) \right)
\end{equation*}%
holds.
\end{corl}

\section{Conclusion}

We have introduced some interesting inequalities, including Chebyshev's
inequality, H\"{o}lder's inequality and Minkowski's inequality for universal
integral on abstract spaces. Furthermore, the reverse previous inequalities
for semiconormed fuzzy integrals are presented. For further investigation,
it would be a challenging problem to determine the conditions under which (%
\ref{eq3-2}) becomes an equality.

\end{document}